\documentclass[10pt, US Letter, conference]{ieeeconf}

\usepackage[english]{babel}
\usepackage{float}
\usepackage{graphicx}
\usepackage{amsmath}
\usepackage{xcolor}
\usepackage{MnSymbol}%
\usepackage{wasysym}%
 \usepackage{amsfonts}

\usepackage{enumerate}

\usepackage{bm}

\newtheorem{thm}{{Theorem}}

\newtheorem{cor}[thm]{Corollary}

\newcommand{\Cmnt}[1]

\newcommand{\TR}[2]{#1}  

\newcommand{\eop}{\hfill{$\blacksquare$}}

\usepackage{float} 
\graphicspath{{pics/}}

 \newtheorem{lemma}{Lemma}
 
 \title{write up}

 \newcommand{\indc}[1]{ {\cal X}_{ \{#1\}}}
\newcommand{\vcmnt}[1]{ #1}
\begin{document}
\title{Asymmetric  Information Acquisition Games   }
 \author{ Vartika Singh and  Veeraruna Kavitha,    \\
 IEOR, Indian Institute of Technology Bombay, India   }
 \maketitle

 \renewcommand{\a }{{a}}
 \newcommand{\dtau}{\delta}

\newcommand{\btau}{ {\bar \tau}}
\newcommand{\ba}{ {\bar a}}
\renewcommand{\a}{{a}} 
\newcommand{\x}{{x}}

\newcommand{\z}{{\bf z}} 
\renewcommand{\t}{{\bm \tau}}

\begin{abstract}

We consider a  stochastic game with partial,  asymmetric and non-classical information, where the agents are trying to acquire as many available opportunities/locks as possible.  Agents have access only to local information, the information updates are asynchronous and our aim is to obtain relevant equilibrium policies. Our approach is to consider optimal open-loop control until the information update,   which allows managing the belief updates in a structured manner. The agents continuously control the rates of their Poisson search clocks to acquire the locks, and they get rewards at every successful acquisition; an acquisition is successful if all the previous stages are successful and if the agent is the first one to complete.   However, none of them have access to the acquisition status of the other agents, leading to an asymmetric information game.  Using standard tools of optimal control theory and Markov decision process (MDP)  we solved a bi-level control problem; every stage of the dynamic programming equation of the MDP  is solved using optimal control tools. We finally reduced the game with an infinite number of states and infinite-dimensional actions to a finite state game with one-dimensional actions.  We provided closed-form expressions for Nash Equilibrium in some special cases and derived asymptotic expressions for some more.

 \Cmnt{
We consider a  stochastic game with partial,  asymmetric and non-classical information.  Agents have access only to local information, the information updates are asynchronous and our aim is to obtain relevant equilibrium policies. Our approach is to consider optimal open-loop control until the information update,   which allows managing the belief updates in a structured manner. The agents continuously control the rates of their Poisson search clocks to acquire the locks, and they get rewards at  every successful acquisition; an acquisition is successful if all the previous stages are successful and the agent is the  first one to complete. They also derive a terminal reward, upon successful  completion of the project. However none of them have access to the acquisition status of the other agents, leading to an asymmetric information game.  Using standard tools of optimal control theory and Markov decision process (MDP)  we solved a two-level control problem; every stage of the dynamic programming  equation of the MDP  is solved using optimal control tools. We finally reduced the game with an infinite number of states and infinite dimensional actions to a finite state game with one-dimensional actions.  We provided closed form expressions for Nash Equilibrium in some special  cases, and derived asymptotic expressions for some more. 
 }

\end{abstract}

\vspace{-3mm}
\section{Introduction}
 
 We consider   non-classical information games (as described in \cite{basar})  inspired by the  full information  games  considered in \cite{eitan}.
In \cite{eitan},   agents attempt to acquire     $M$ available destinations;  each agent controls its rate of acquisition (advertisement through a Social network)  to increase its chances of winning the destinations, while trading off the cost for acquisition.   They considered full information games, wherein  the agents at any time  know the number of destinations  available (not  acquired by some agent);  they also  considered   no-information games, where the agents would not even   know the number of destinations/locks acquired by themselves. 
 
 It is more realistic to assume that the agents  have access only  to  local information;   they  might  know the number of locks acquired by themselves, but  not the numbers acquired by  others.  This leads to partial, asymmetric and non-classical information games;    Basar et. al in \cite{basar} describe a game to be of non-classical information type, and we describe the same in our   words:  if  
 the state of  agent $i$  depends upon the actions of   agent $j$, and if agent $j$  has  some information which is not available to agent $i$ we have a non-classical information game.
 These kind of games are typically hard to solve (\cite{basar}); when one attempts to find best response against a strategy profile of others, one would require belief of  others states, belief about the belief of others, and so on.

 We have some initial results related to this asymmetric information game in \cite{Mayank}, where we solved the problem for the case with two destinations  and  two agents.  We also considered that the 
 locks (represent destinations) have to be acquired in a given order and only the  agent that  first wins all the locks would receive a (single) prize.  We now consider a significant generalization of the   game:  a) we consider general number of locks and agents; b) the agent that wins all the first $k$ locks before the others would receive prize $c_k$ and this is true for any lock. 
 
Our approach to this problem can be summarized as ``open loop control till information update" (as in \cite{Mayank}).  With no-information, one has to resort to open loop policies (action changes with time, but, oblivious to the state). This is the best when one has no access to information updates.  With full information one can have closed loop policies (actions can depend on state).  Further, in full information    controlled  Markov jump processes, every agent is informed immediately of the jump in the state and can change its action based on the change.  In our case we have access to partial information,  the agents can observe only some jumps and not all; 
thus
we need policies that are open loop type till  an information update.   At every information update, one can choose a new open loop control depending upon the new information.

The agents have no access to the information   of  others,  however upon contacting a lock they would know if they are the first to contact.
Any strategy profile  in our  game is described by one open loop policy for each state, each state is primarily described by the time of acquisition of the previous lock; thus    we have   an infinite dimensional (strategy space is polish space) game. 
We used the tools of optimal control theory (Hamiltonian Jacobi equations) and   Markov decision process to reduce this  infinite dimensional game to a one dimensional game; we showed that  the best response against any strategy profile of others includes a time threshold policy (acquisition attempts to be made with full intensity till a threshold of time, after which  there would be no attempts); more importantly we showed that these thresholds can be specified by a single constant (for each lock), irrespective of the time of acquisition of the previous lock. We finally showed that the {\it reduced game is a strict concave game, proved the existence of unique Nash equilibrium}   and provided a simple algorithm to compute the same.  
Our results matched with those in \cite{Mayank}, for all the cases considered   there. We   have closed form expressions (some are asymptotic) for the Nash equilibrium for some cases. 

\Cmnt{
We considered one and two  lock acquisition problems, any agent wins reward one if it acquires all the locks and if it is the first one to acquire the locks.  The agents have no access to the information/state  of the others,  however upon contacting a lock they would know if they are the first to contact.  We obtained Nash equilibrium for these partial, asymmetric information games; a pair of (partial) state-dependent time threshold policies form Nash equilibrium. We obtained these results (basically best responses) by solving Dynamic programming equations applicable to (partial) information update epochs  and each stage of the Dynamic programming equations are solved by solving appropriate optimal control problems and the corresponding Hamiltonian Jacobi equations. 
 
A block chain network is a distributed ledger that is completely open to  all nodes; each node will have a copy of transactions (in case of currency exchange).    If a new  transaction is  to be added (linked) to the previously existing chain in the form of a new block, it  requires the miners (designated special nodes) to search for a key (encryption), that enables it to be added to the ledger. This search of the key requires computational power and time. The first miner to get the right key,   gets a financial reward.  If  the miners would not   know the status of the search efforts of others, the resulting game is exactly as in our one lock problem. Two lock problem can be viewed as the extension of one lock problem, wherein a second key is required to gain the reward.  
 }

Our models can capture a variety of applications, e.g., social network problem as in \cite{eitan} with one destination or a block chain problem as in \cite{Mayank}.  For general $M$, we solve problem of winning a project  (unaware of others success) with multiple completion phases. 

\section{Problem statement}
 There are $n$ agents   competing to win a project. The project is to acquire $M$ locks before the deadline $T$. The locks are ordered, i.e, all the agents will compete  for the first lock  in the beginning, after  which  they will compete for the second lock and so on. The agent that contacts all the locks successfully wins the project and gets a terminal reward. Further  we say $k$-th  contact (of $k$-th lock) is successful if this contact
  happens before time $T$ and if  the  agent was  the first one to contact all the previous ($k-1$) locks and the $k$-th lock;  the agent  gets some   reward at every successful contact. The acquisition/contact process of each agent is modelled by independent  (possibly non-homogeneous) Poisson Processes;  they  can choose the rate of their (Poisson) contact process    as a function of time $t \in [0,T]$.
  A higher rate of contact   will increase  the chances of success but will also incur higher cost. The aim of the agent is to maximize its expected reward. 

 \noindent {\bf Information structure: } The agents have partial/asymmetric information about  the number of locks acquired by various agents and would use the available information to design their rate functions optimally. Any agent  would know at all the times information related to  its contact attempts, however has limited access to that of  the others. When  it  contacts a  lock, it would know if it is successful;  if not, it gets to know  that  it  is not the first one to contact. So, agent gets some partial information about the state of others.

\noindent{\bf Decision epochs: } As already mentioned, an agent has access only to partial information. 
There is a  (major) change     in the available information at  the lock-contact epochs\footnote{By convention, the start of the process commences with 0-th contact.};  it would know 
successful/unsuccessful status of  the   lock  immediately after contact, and based on this information  the agent can choose an optimal action. Hence contact epochs form the natural decision epochs. 
 Further, these epochs will be exponentially distributed random variables with the parameters chosen by the agents, so it is clear that the  decision epochs of different agents will not be synchronous. \par

{\bf State:} 
The state  of  any  agent at  any  time $t\in[0,T]$ is the information available to it at that  time. 
 The information available to the agent $i$ after contacting $(k-1)$-th lock,  $z_k^i$, has two components:  i)  a flag $l^i_k$ represents whether the contact was successful ($l^i_k=1$) or unsuccessful ($l^i_k = 0$); and ii) the time of $(k-1)$-th contact denoted by $\tau^i_{k-1}$. Thus $z_k^i=(l^i_k,\tau^i_{k-1})$   is the state of  agent $i$ at decision epoch $k$. The state   remains the same in the time interval $[\tau^i_{k-1},\tau^i_k)$. The initial state $z_1^i\vcmnt{=(l^i_1,\tau^i_0)}$   is simply  set to $(1,0)$;  the process  starts at 0 and $l_1^i  $ is set as 1 for convenience. 

{\bf Actions: } 
The agents choose the rate functions (defined till $T$) at their respective decision epochs based on their states; these functions  are \textit{open loop policies}, wherein 
the dynamic action  is independent of the state of the system; the agents change their rate function only  at  next decision epoch. Such an approach is called ``Open loop policy till information update" as in \cite{Mayank}. 
 The rate of contact, for agent $i$,  at any time  can take values in the interval $[0, \beta^i]$ and  the rate function   is  measurable.   To be precise, agent $i$ at decision epoch $k$, i.e., at time instance  $\tau^i_{k-1}$,  chooses  an action  ${a}_{k}^i \in L^\infty [\tau^i_{k-1}, T]$, as the acceleration process to be used till the next acquisition. Here  $L^\infty [\tau^i_{k-1}, T]$ is the space of all measurable functions that are  uniformly bounded by the given constant;   the bounds ($\beta^i$ for agent $i$) can be different for different agents, nevertheless we avoid notation $i$ for simpler representation.   These form a closed subset of  Polish space\footnote{Complete and separable space.} of 
 essentially bounded functions, i.e.,  the space of functions with   finite essential supremum norm: 
 $$
|| \ a ||_\infty :=  \inf \{ \beta:  |a(t)| \le \beta  \mbox{ for almost all } t \in [\tau^i_{k-1}, T]\}.
 $$

 {\bf Strategy: }The strategy of player $i$ is a collection of open loop policies, one for each state and lock,  as given below:
 
 \vspace{-4mm}
 {\small
\begin{eqnarray}
\label{Eqn_Strategy}
 \pi^i =  \left  \{  \ a_{k}^i (\cdot \ ;\ z^i_k) \in L^\infty ; \mbox{ for all  }  z^i_k  \mbox{ and all } {k \in \{1, \cdots, M\} } \right  \},
\end{eqnarray}}where $\ a_1^i ( \cdot  ;\ z^i_1)$ represents the open loop policy/rate function used at  start,  while $\ a_k^i ( \cdot \ ;\ z^i_k)$ represents the  same to be used after  $(k-1)$-th contact;  this  choice depends upon   the  available information $ z^i_k.$
At times, notation $z^i_k$ is dropped   and we use $L^{\infty}$ in place of $L^\infty [\tau^i_{k-1}, T]$, to simplify the  explanations.  

Our work majorly  analyses   best response (BR) of players, hence we introduce the following notations. Without loss of generality  we consider BR of  agent $i$. 
{\it Let $N:=\{1, 2, \cdots, n\}$ be the set of players and $\pi^j$ be the policies of  all other players represented by     $j:=-i  := N - \{i\} $.} 

{\bf Rewards/Costs: } The reward of  agent $i$ is $c^i_k$, if it contacts all the first $k$ locks successfully   (before time $T$ and before other agents), while the terminal reward is $c_M^i$.   This implies,  an agent   unsuccessful with first contact,   has  no incentive to attempt for  further locks; so it would  remain silent henceforth. 
 
Recall that     $ T\wedge \tau_k^i $  represents the time  instance\footnote{ The contact clocks $\{\tau^i_k\}$ are free running Poisson clocks, however we would be interested only in those contacts that occurred  before deadline $T$.} till which the  $k$-th  lock is attempted.   Then the cost spent on acceleration  for the  $ k$-th contact equals,  
\begin{eqnarray}
\bar{a}_k^i (T\wedge \tau_k^i ), \mbox{ with }   \bar{a}^i_k  (t)   :=  \int_{\tau^i_{k-1}}^{t} a^i_k ( s )  ds. \label{Eqn_bar_a}
\end{eqnarray}
 Thus the expected (immediate) utility for stage $k$  equals:
 \begin{eqnarray} 
r_k^i(z_k^i,  a_k^i; \pi^j) =  \left \{ \begin{array}{lll}
c^i_k  P_k^i(z^i_{k}; a_k^i; \pi^j)  &   \mbox{if } l_k^i = 1 \mbox{ and}  \hspace{-4mm}\\ 
\hspace{1mm} - \nu  E\left [ \bar{a}_k^i (T\wedge \tau_k^i ) | z_k^i \right ] & \ \ \ \tau_{k-1}^i  < T ,  \hspace{-4mm}  \\
0    & \mbox{else, } 
\end{array}   \right . \hspace{-1mm} \label{Eqn_run_cost}
\end{eqnarray}where $P_k^i $   represents 
the probability of   successfully contacting $k$-th lock, conditioned on $z^i_{k} = (1, \tau^i_{k-1})$, i.e., conditioned that $(k-1)$-th lock is   acquired successfully at $\tau^i_{k-1}$;  and $\nu > 0$ is the trade-off factor between the reward and the cost. 

When $k=1$ the probability of success $P_k^i$ also depends upon the failure of other agents, i.e., depends upon  $\pi^j := \pi^{-i}$:  
\begin{eqnarray*} \hspace{2mm}
P_1^i (z^i_{1}, a^i_1; \pi^j)  =  \int_0^T   \eta_1^i (s; \pi^j)    e^{ - {\bar a}_1^i (s) }   a_1^i (s)  ds, \end{eqnarray*}where the probability of failure of other agents before time $t$ equals ($\mathcal{X}$ is the indicator and see \eqref{Eqn_bar_a} for definition of $\bar{ a}$):
\begin{eqnarray}\label{eqn_failure}
\eta_k^i (t; \pi^j) =  \indc{k > 1}  + \indc{k=1}  e^{-\sum_{m \ne i} {\bar a}_1^m (t) }. 
\end{eqnarray}
In the above, the indicators are introduced to have unified notation;   
for $k > 1$, the probability of  success $P_k^i$  conditioned on success till $(k-1)$-th lock (now there is no opposition),  equals: \vspace{-4mm}
\begin{eqnarray}\label{eqn_eventual_success}\hspace{10mm}
P_{k}^i (z^i_{k}, a_k^i; \pi^j) = \int_{\tau^i_{k-1}}^T\eta_k^i  (s;\pi^j)  e^{ - {\bar a}_k^i (s) }   a_k^i (s)  ds.
\end{eqnarray}

It is easy to observe that for any given $\ a^i_k (\cdot) \in L^\infty$, the expected cost equals  (see (\ref{Eqn_bar_a}) and with $\tau_0^i := 0$):

\vspace{-4mm}
{\small \begin{eqnarray}
\label{Eqn_Cost}
\vcmnt{E[\bar{a}_k^i ( T \wedge \tau^i_{k}) | z_k^i  ] } = 
 \bar{a}_k^i (T ) e^{ - \bar{a}_k^i (T ) }  
 + \hspace{-2mm}  \int^{T }_{\tau_{k-1}^i}\hspace{-2mm}   \bar{a}_k^i (s)      e^{   - \bar{a}^i_k  (s)  }  a^i_k ( s)  ds.   
\end{eqnarray}}If  an agent  fails to contact   before the deadline $T$, it still  has to pay for the entire duration $T-\tau_{k-1}^i$  and hence the first term in the above equation.

{\bf Game Formulation: } This problem can be modelled as a strategic/normal form non-cooperative game, 
$${G}=\langle  {N}, {S}, {\Phi} \rangle, \mbox{ with } {S} = \{S^i\}_i \mbox{ and } {\Phi} = \{\phi^i\} ,$$ where 
the set of players
$N=\{1,2,\dots,n\}, $
the strategy set of player $i$ is the class of  all possible  strategies as in (\ref{Eqn_Strategy}), i.e., $S^i = \{\pi^i\}$, and the (overall) utility  of agent $i$   is given by
\begin{eqnarray}
\label{Eqn_Objectivei}
\phi^i (\pi^i, \pi^j ) =  \sum_{k=1}^M  E[ r_k^i(\ z_k^i, \ a_k^i;  \  \pi^j  )  | z^i_1 = (1,0) ].
\end{eqnarray}
Our aim is to find a  tuple of strategies (that depend only upon the available  information)  that form the Nash equilibrium (NE). 
We conclude this section by mentioning some of the main results of the paper. 

 \subsection{Important results}
 \label{sec_Important_results}
 
 We solved this problem for $M$ locks and $n$ players;     the problem is converted  to a much simplified and   reduced strategic form game,  such that   the (unique) Nash equilibrium of  the reduced game  is  also the   Nash equilibrium of the original game.   We found the reduced game with the help of following results.  Before stating  the  results, we  require the following definitions. 

\textit{Threshold Policy:} 
A threshold policy is an open loop policy, which takes value $\beta$ (maximum possible value) in the  interval $[s,\theta]$ (if $s\le \theta$) and value zero in the remaining interval, where $s$ is the starting point of the control and $0 \leq \theta \leq T$.  
Threshold policies are  represented by  $\Gamma_{\theta; s}$, where (for any $s \le t \le T$),  \vspace{-6mm}
\[ \hspace{10mm} \Gamma_{\theta ; s}  (t) = \begin{cases} 
      \beta & \mbox{if }  t\leq \theta, \\
      0 & \mbox{else.}
   \end{cases}
\]
The  rate function is 0 for all $t$, if the starting point $s > \theta$.

{\it Threshold (T) strategy:} is a strategy   made up of threshold policies. A  typical  T-strategy is defined by threshold functions $\{ \theta_k ( \cdot  ) \}_k$ and is defined as below:

\vspace{-4mm}
{\small \begin{eqnarray*}
\pi \hspace{-2mm}&=&\hspace{-2mm} \{ \theta_1 ( \cdot )  \dots, \theta_M ( \cdot ) \},   \\
\hspace{-2mm} & = &  \hspace{-2mm}  \big   \{  a_{k} (\cdot \ ; z_k) ;  a_k (\cdot\ ; z_k) =  \Gamma_{\theta_k (s)  ; s } (\cdot)      \mbox{ when  }  z_k = (1, s)  \mbox{ and }    \\
&& \hspace{13mm}       a_k (\cdot\ ; z_k) =  \Gamma_{0; s} \mbox{ when } z_k = (0, s) \mbox{, for any }  \ k  \big \}.
\end{eqnarray*}}Basically when the state of player $i$ is $z_k^i = (l^i_k, \tau_{k-1}^i)$ with $\tau^i_{k-1} = s$  and $l^i_k = 1$, then starting from time  $s$   player $i$ uses  threshold policy with threshold equal to  $\theta_k (s)$  as the open loop rate function;    when $l^i_k = 0$,  we equate  the thresholds   to 0, i.e., the player stops trying any further.  
Many a times the {\it starting point is obvious and hence the corresponding notation is dropped from the subscript.} \\
{\it M-Thresholds  (MT) Strategy:}  This is a special type of    T-strategy in which the  thresholds for any given lock remain the same irrespective of the time of acquisition   of the previous lock. 
A  typical  MT-strategy is defined by  $M$-thresholds: 

\vspace{-4mm}
{\small \begin{eqnarray*}
\pi &=&  \{ \theta_1 \dots, \theta_M \} ,
 \\
 & = &  \{   a_{k} (\ \cdot \ ;\ z_k)   ;  \  a_k (\ \cdot\ ; z_k) =  \Gamma_{\theta_k}  (\cdot)      \mbox{ for any  }  z_k  \mbox{ and }  k    \}.
\end{eqnarray*}}
In other words, MT-strategy is  completely represented by $M$-thresholds $\{ \theta_k \}_{k\le M}$,  one for each $k$: 
a) here  $\theta_k$  represents  the time threshold till which agent can attempt to acquire  $k-$th lock; b)   threshold  $\theta_k$ is independent of  $\tau_{k-1}$, the time of acquisition of the $(k-1)$-th lock; and  c) if   $\tau_{k-1}$ is bigger than $\theta_k$ (i.e., if the $(k-1)$-th lock is acquired after the threshold for the $k$-th lock) then the agent would no longer attempt to acquire the $k$-th lock under the given MT-strategy.

To summarize,  the strategy of any player is made up of open loop policies,   one  for every (possible) state and for each decision epoch.  Thus we have infinite dimensional actions, however  
  the following structural results about the best response strategies reduce the complexity of the game significantly:

\begin{thm}
{\bf [Threshold Strategy]}
\label{Thm_threshold_policy}
There exists 
a (threshold) $T$-strategy that is a  best response  strategy  of any given  player, against any given  strategy profile  of opponents.   \eop
\end{thm}{}

\begin{thm}{\bf[$M$-thresholds strategy]}
\label{thm_monotonocity_of_thrshold}
There exists    an MT-strategy  that is a  best response  strategy  of any given  player, against any given  strategy profile  of opponents.   \eop
\end{thm}{}

The proofs  of  both the Theorems are in the next section,  section \ref{sec_proof_of_Theorems}.
By virtue   of the above   Theorems,
 there exists  a best response   strategy   represented completely  using $M$ thresholds, which is optimal among (uncountable)   state dependent strategies with   infinite dimensional  strategy   space.     
 We will thus have a reduced  game  in ${\cal R}^M$ which is further analyzed in section \ref{sec_ReducedGame}; we find NE of the original game among  these  MT  strategies.  The unique NE of the reduced game is characterized by Theorems \ref{Thm_Upsilon_two} and \ref{thm_unique_NE} of section \ref{sec_ReducedGame}.

\section{Best responses}
\label{sec_proof_of_Theorems}

\Cmnt{{\color{red}
Emphasize BR against $\pi^j$

Show that every stage the DP equations can be modelled as optimal control problem (number it)

The HJB PDEs are common tools to derive solutions of these kind of optimal control problems..
Write the corresponding HJB PDE  (cite flemming) ..

Theorem 4: for every state, and every $z_k$, i) there exists a   unique viscosity solution which is continuous for PDE (..)  
ii) There exists a policy that optimizes the  optimal control problem. End of Theorem

So, every stage of DP equation has an optimizer (when BR is being is computed) thus by Theorem 4.5.1. of Puterman the policy $\pi_*^i = \{a^{i*}_k (z_k) \}$ forms the BR against $\pi^j$.
}
}
Our aim is to derive Nash equilibrium (NE) for this partial and asymmetric  information stochastic game.
We begin with deriving the best response (BR) of player $i$ against  any given  strategy profile $\pi^j$ of the opponents.

\noindent {\bf Dynamic programming equations}
The BR
is obtained by maximizing  the objective function
 (\ref{Eqn_Objectivei}) with respect to the strategies  $ \pi^i \in S^i$.  It is easy to observe that this  optimization is an example of a Markov decision process 
which 
 can be solved using ($M$-stage) dynamic programming (DP) equations given below  (see \cite[Theorem 4.5.1 and Remarks after]{Puterman} which covers the results for Polish spaces): 

\vspace{-4mm}
{\small \begin{eqnarray}
\nonumber 
v^i_{k} (z^i_k; \pi^j) &=&\hspace{-3mm} 0  \hspace{3mm}\mbox{ if }  k > M \mbox{ or if } \tau_{k-1}^i >  T \mbox{ or if }  z^i_k=(0, \tau^i_{k-1}), \label{eqn_dp_2}  
\end{eqnarray}{\normalsize and other wise  the value function equals,}
 \begin{eqnarray}
v^i_k (z^i_k; \pi^j) =\sup_{\ a_k^i \in L^\infty } \left \{ r_k^i (z^i_k,a^i_k; \pi^{j})+ 
  E[  v^i_{k+1} (z^i_{k+1}; \pi^{j}) | z^i_k, a^i_k] \right \}.    \hspace{-8mm} &
     \nonumber \\
&  \label{Eqn_DP}
\end{eqnarray}}
Observe that (\cite{Puterman} and see equation (\ref{Eqn_Objectivei}))
$$
v_1^i (z_1^i ; \pi^j) =   \sup_{ \pi^i \in S^i }  \phi^i (\pi^i;  \pi^j ),
$$and thus the BR is obtained by solving the DP equations. 
The $k$-th stage DP equation can be re-written as below:
\begin{equation}\label{eqn_value_func}
    v_k^i(z^i_k; \pi^j)  \ =  \sup_{\ a_k^i \in L^\infty } J^i_k ({z^i_k}, a_k^i; \pi^j) ,  
\end{equation}
where the cost $J^i_k$ is defined by (see equations \eqref{Eqn_run_cost}-\eqref{Eqn_Cost}):

\vspace{-4mm}
{\small
\begin{eqnarray}\label{cost}
    J^i_k ({z^i_k}, a_k^i; \pi^j)  \hspace{-3mm} &=&\hspace{-4mm}     \int_{\tau^i_{k-1}}^{T}   \left  ( h_k^i (t) -\nu \bar{a}_k^{i}(t) \right  ) 
  a_k^i ( t) e^{-\bar{a}_k^{i}(t)} dt  \nonumber \\ && \hspace{28mm}
 -\nu \bar{a}_k^{i}(T) e^{-\bar{a}_k^{i}(T )},  \mbox{ with } \nonumber \\
 h_k^i (s)\hspace{-3mm}  &:=& \hspace{-3mm}\big (c^i_k + v_{k+1}^i( (1,s); \pi^j) \big) \eta_k^i(s;\pi^j).
\end{eqnarray}}
\noindent{\bf Optimal control:}  
From the structure of the optimization problem  \eqref{eqn_value_func} defining the $k$-th stage DP equation, it is clear that one can solve it using an appropriate optimal control problem. 
One can write this  optimization problem as:
\begin{eqnarray*}
v_k^i(z^i_k; \pi^{j})  &=&  u(\tau, 0) ; \mbox{ when } z^i_k = (1, \tau) ,
\end{eqnarray*}
where   $u(s, x)$ (defined for any $s \in [\tau, T]$  and any $x$) is  the value function
of 
the  optimal control problem  with details: 

\vspace{-4mm}
{\small\begin{eqnarray}
u(s, x) \hspace{-2mm} &:= & \hspace{-5mm} \sup_{a \in L^\infty [s, T] } \hspace{-2mm} J (s, x,  a ),   \mbox{ where the objective function, }   \hspace{8mm}\label{eqn_opt_ctrl}  \\
 J   (s, x,  a ) \hspace{-2mm} & :=& \hspace{-3mm}\int_{s}^{T} \left ( h_k^i  (s')    - \nu   x (s')  \right ) a ( s') e^{-x(s')} ds'   +g(x(T )),  \nonumber 
\end{eqnarray}}with state process of  the optimal control problem given by (for any~$s' \in [s, T]$)
$$
\stackrel{\bm\cdot}{x}  (s') =    a  ( s'), \mbox{ with } x(s)= x  \mbox{ and thus }  x(s') = x+ \int_s^{s'}\hspace{-2mm} a ( {\tilde s}) d {\tilde s} 
$$ and
 with terminal cost  \vspace{-6.5mm}
\begin{equation*}
\hspace{5mm}
    g(x) = - \nu x e^{ -x }. \nonumber
\end{equation*}

We need to solve this optimal control problem to get the BR, 
and the standard technique to solve such  problems is using  Hamiltonian Jacobi  (HJB) PDEs \cite{feller}, and the one corresponding to   (\ref{eqn_opt_ctrl}) is given  by:
 
 \vspace{-4mm}
 {\small\begin{eqnarray}
u_s  (s, x) + \sup_{a \in [0, \beta^i]}   \ a\{  ( h_k^i (s)    - \nu   x   )e^{-x}    +  u_x (s,x) \}   \hspace{-1mm}&  \hspace{-1mm}= \hspace{-1mm}&  \hspace{-1mm} 0, \mbox{ with }  \nonumber   \\
  u(T, x)     \hspace{-1mm}&  \hspace{-1mm}= \hspace{-1mm}&  \hspace{-1mm}   g( x),  \hspace{7mm}\label{eqn_hjb_pde}
\end{eqnarray}}where $u_s$, $u_x$ are partial derivatives of the (optimal control) value function. Using the standard tools we immediately obtain the following result (Proof in Appendix):
\begin{thm}{\bf [Existence]}
\label{thm_existance}
At any stage $k$ and  any   state  $z^i_k$, 
\begin{enumerate}[(i)]
\item The  optimal-control value function $u(\cdot, \cdot)$ is the  unique viscosity solution of the HJB   PDE \eqref{eqn_hjb_pde}.
Further, the value function $u (\cdot, \cdot)$  is   Lipschitz continuous in $(s,x)$.
\item We have a $a^* (\cdot) $ that  solves  the    control problem.  \eop
\end{enumerate}
 \end{thm}

{\bf Remark:} Theorem \ref{thm_existance} implies every stage of DP equation has an optimizer, i.e., for every lock $k$ and state $z^i_k$, the player $i$ has a best response policy, call it 
$a_k^{i*}(\ \cdot \ ; z^i_k)$. Thus by  \cite[Theorem 4.5.1 and the following remarks about Polish space]{Puterman}, 
$\pi^{i*} =\{ a_k^{i*}  (\ \cdot \ ;z^i_k)\}_{k, z^i_k }$ is a  BR  strategy of player $i$ against $\pi^j.$  

We further show that the optimal strategy  (in BR)  can be  a threshold strategy in the following. We begin with:
\begin{lemma}\label{lemma_earlier_control}
For any  $k$,  $z^i_k$ and open loop policy $a^i_k(\cdot)$,  one can construct a threshold policy $\Gamma_{\theta } (\cdot)$ such that
 the contact time under threshold policy (denoted by  $\tau_\theta  $) stochastically dominates   that ($\tau_a  $) under $a^i_k(\cdot)$, 
  i.e,  
$
\tau_\theta   \stackrel{d} {\le }   \tau_a,$  i.e., for any monotone  decreasing function $f$, 
$$  E[f(\tau_a  ) ]  \le E[f(\tau_{_\theta}   ) ] .
$$ Further  the expected costs (\ref{Eqn_Cost}) under both the policies are  equal.  Also, the probability of successful contact $ P_k^i$   given in \eqref{eqn_eventual_success} is better with  the threshold policy, i.e.,
$$  P_k^i(z^i_{k};a^i_k(\cdot); \pi^j)  \le  P_k^i(z^i_{k}; \Gamma_\theta; \pi^j)   . $$
 
\end{lemma}{}
{\bf Proof} is in Appendix.  \eop 

Thus  at any stage $k$  and  for any state $z^i_{k}$, agent $i$ can  contact the locks  faster using threshold policies and the running costs  \eqref{Eqn_run_cost} are better.
When the locks are contacted faster, i.e., when the next stage starts earlier,  the ``cost to go" (the value function) till the end is better   as shown below:

\begin{lemma}
\label{Lemma_value_function_mono} 
For any lock $k$, the value function with $t \le \tau$, 
$$
v^i_k (z^i_k;\pi^j)  \ge v^i_k ( {\bar z}^{i}_k;\pi^j)  \mbox{ when }  z^i_k = (1, t) \mbox{  and } {\bar z}^{i}_k = (1,  \tau) . 
$$  
\end{lemma}{}
{\bf Proof} is in Appendix.  \eop 
 
 \subsection{Completing the proof of Theorem \ref{Thm_threshold_policy} }

 We need to prove that for every lock,  there exists a threshold policy which is in best response against any fixed strategy of the opponents. From Theorem \ref{thm_existance} we have the existence of the optimal policy (BR policy against any fixed $\pi^j$) for every lock. Let's denote the optimal policy for acquiring $k$-th lock by $a^*(\cdot)$, and, say it is not a threshold policy. Let $\tau_{a^*}$ denote the contact epoch under policy $a^*$. If $a^*$ is the optimal policy, it must be the optimizer (maximizer) of the  $k$-th stage DP equation (\ref{eqn_dp_2}) and so,
  $$v^i_k (z^i_k; \pi^j) = 
 r_k^i (z^i_k,a^*; \pi^{j}) + E[  v^i_{k+1} ((l,\tau_{a^*}); \pi^{j})  | z^i_k,  a^*] .$$
Construct a threshold  policy $\Gamma_\theta$ using $a^*$ as in proof of Lemma \ref{lemma_earlier_control} and then from the same Lemma:
$$ r_k^i (z^i_k,a^*; \pi^{j}) \le  r_k^i (z^i_k, \Gamma_\theta; \pi^{j}) ;$$
from Lemma \ref{Lemma_value_function_mono}  the value function is a  non-increasing function of time, and hence again using Lemma~\ref{lemma_earlier_control} we have,
$$ E[  v^i_{k+1} ((l,\tau_{a^*}); \pi^{j})  | z^i_k,  a^*]\le  E[  v^i_{k+1} ((l,\tau_{\theta } ); \pi^{j})  | z^i_k,   \Gamma_\theta].$$
This proves  that  a threshold policy is among BR policies against $\pi^j$, using \eqref{Eqn_DP}. \eop
\Cmnt{
 $$\bar{v}^i_k ( z_k; \pi^{j}) = 
 r_k^i ( z_k,  a; \pi^j) + E[  v^i_{k+1} ( (l,\tau_a); \pi^j) ) |  z_k,  a] $$
Note that the first term of the DP equation is the instantaneous reward, for $k>1$ the instantaneous reward will be same for both the policies $a$ and $a'$. For $k=1$, it follows from part 2 of lemma \ref{lemma_earlier_control}, the instantaneous reward will be more  in the threshold policy 
 \begin{equation}\label{eqn_rew_theshold}
     r_k^i (\ z_k, a; \pi^j) \geq r_k^i (\ z_k, \ a'; \pi^j)
 \end{equation}{}

if $X = \tau_a $ is stochastically bigger {\color{red}smaller} than $Y = \tau_{a'}$, i.e., if  $X \stackrel{d}{\ge} Y$ and $f = v_{k+1} ( (l, . ), \pi^j)$  is any {\color{red} non-decreasing ?? or non-increasing} function then 
$$
E[f(X)] \ge E[f(Y)]. 
$$
{\color{red}
$$
E[f(X)] \le E[f(Y)]. 
$$
}
 
By using part 1 of lemma \ref{lemma_earlier_control}, the time of contact of the $k$th lock in the threshold policy will be less,i.e., the kth lock will be contacted earlier using threshold policy with probability 1,
 $$P(\tau_a\leq \tau_{a'})=1 \mbox{ and }P(\tau_{a}\ < \tau_{a'})>0
 $$
 so  by the lemma \ref{Lemma_value_function_mono},
 $${ v}^i_{k+1} ( (l,\tau_{a'}); \pi^j) )\leq \bar{ v}^i_{k+1} ((l,\tau_{a}); \pi^j) )\ \ \ \  \ \ \ \ a.s.$$
{\small\begin{equation}\label{eqn_value_theshold}
    E[  v^i_{k+1} ((l,\tau_{a'}); \pi^j) ) | \ z_k, \ a_k]\leq E[  v^i_{k+1} ((l,\tau_{a}); \pi^j) ) | \ z_k, \ a'_k] \ a.s.
\end{equation}}
from equation (\ref{eqn_rew_theshold}) and equation (\ref{eqn_value_theshold}), we have
$$\bar{v}^i_k ( z_k; \pi^j) \geq v^i_k (z_k; \pi^j)$$
 this proves that the policy $a$ which is a threshold policy is optimal.   \eop

 Thus the above Theorems confirm    there exists a BR policy  which is specified by  one threshold for each $k$ and is sufficient to obtain $a_k^{i*} ( ; z_k) $ when $z_k = (1, 0)$, against any strategy of opponents $\pi^j$. This implies the existence   
Using the above two results we conclude the  proof of our first main result, Theorem \ref{Thm_threshold_policy}.}

 \subsection{Proof of Theorem \ref{thm_monotonocity_of_thrshold}}
 We prove the second Theorem in two steps;  in the first step we show that the optimal policies corresponding to search of  $k$-th lock  coincide in all possible time intervals, irrespective of the start of this search, $\tau^i_{k-1}$, in the following sense:
 \begin{thm}
\label{theoerem_coincidence}
Let $\tau\geq t$. 
The optimal/BR policy to acquire $k$-th lock, $a^{i*}_{k}(\cdot\ ; z )$ with $z = (1,t)$ coincides with  BR policy $a^{i*}_{k}(\cdot\ ;  z')$  with $z' = (1, \tau)$ from $\tau$ onwards,   i.e.,
$$\hspace{20mm}a^{i*}_{k}( s; z)=a^{i*}_{k}( s; z') \mbox{ for all }  \tau  \le s \le T. \hspace{20mm} \TR{}{\mbox{ \hspace{-15mm}\eop}} $$
\end{thm}
\TR{
{\bf Proof: } 
The optimal policy to acquire $k$-th lock, $a^{i*}_{k}(\cdot\ ; z )$ with $z = (1,t)$ is optimizer of the value function given in equation (\ref{eqn_opt_ctrl}), i.e., \vspace{-4mm}
$$\hspace{16mm}  u(t,0)=\sup_{a\in L^\infty [t,T]} J(t,0,a).$$
 From  Dynamic programming principle of optimal control problems  \cite[Theorem 5.12]{feller}, we have:
\vspace{-1mm}
{\small \begin{eqnarray}
 \label{Eqn_opt_concatinate}
u(t,0 )   &=   & \sup_{ a \in L^\infty [t, \tau] }  \bigg \{  \int_{t}^{\tau}   \big( h_k^i (s)- \nu   x (s)  \big )  a ( s)   e^{-x(s)}ds  \nonumber   \\
 &&   \hspace{16mm} + u(\tau,x(\tau))  \bigg \} \mbox{ for any }   t\le \tau \le T.  \hspace{4mm}
\end{eqnarray}}
As in \cite[Lemma 4.2]{feller} one needs to find the optimizer for the time interval $[t,\tau)$, considering that the optimal control from $\tau$ onwards will be the same as the one that obtains, the optimal  $ u(\tau,x(\tau))$, where  $x(\tau)$ is the state at $\tau$. 
And if both the problems have optimal policy (the existence for our case is established as in Theorem \ref{thm_existance}), then the optimal policy for the entire interval is given by (as in \cite[Page 10]{feller}):
\begin{eqnarray}
a_k^{i*}(s) =  \left \{ 
\begin{array}{ccccc}
a_1^* (s)  & \mbox{ for all }  s < \tau, \\
a_2^* (s)  & \mbox{ for all }  s > \tau,
\end{array}
 \right . \label{Eqn_policy_stich}
\end{eqnarray}
where $a_2^* (s)$ is the optimal policy  attaining $ u(\tau,x^*(\tau))$,  $x^*(\tau)$ is state at $\tau$  when $a_1^*$ is used in interval $[t, \tau]$ and where $a_1^*$ is the optimizer of \eqref{Eqn_opt_concatinate}.
The optimal control from $\tau$ onwards  in general depend on    state $x(\tau)$ at time $\tau$, but in our case, by Lemma \ref{lemma 3} given in Appendix,  \eqref{Eqn_opt_concatinate} modifies to:
\begin{eqnarray*}
u(t,0)&=&\sup_{ a_k^i \in L^\infty [t, \tau] }\bigg\{ \int_{t}^{\tau} \big( h_k^i (s)- \nu   x (s)  \big )  a_k^i ( s)  e^{-x(s)}ds  \\
&& \hspace{17mm} + e^{-x(\tau)}[u(\tau,0)-\nu x(\tau)] \bigg\}.
\end{eqnarray*}
By   Lemma \ref{lemma 3},    the optimal control from $\tau$ onwards is independent of the state at $\tau$, i.e.,   the optimal control  policies  defining   $ u(\tau,x(\tau))$ and  $ u(\tau,0)$ are the same, and the common one forms a part of $a_k^{i*}$ (see \eqref{Eqn_policy_stich}); this  completes  the proof.\eop

Thus it suffices to optimize for every lock with $z_k = (1, 0)$ (i.e., with $\tau_k = 0$) and the rest of the optimal policies (with different starting time instances) can be constructed using this 
{\it zero-starting optimal policies}, which 
  immediately leads to the following corollary: 
\begin{cor}\label{cor_finite_collection}
The optimal (BR) strategy can be completely specified  by a finite    ($M$)  collection of  control policies,  $ \pi^{i*}_0(\pi^j) := \{a^{i*}_{0_k} (\cdot ) \} $,  one for each lock and each of them  starting at time zero and spanning till time $T$, and such that:
\begin{eqnarray*}
\pi^{i*} (\pi^j)  &=&  \{  \ a^{i*}_{k} (\cdot\  ; z_k ) , \  \mbox{ for all } z_k \},  \mbox{ where } \\
  a^{i*}_{1} (s   ; z_1  )  &=&  a^{i*}_{0_1}  (s),   \mbox { and,  for any }  k > 1 \\
 a^{i*}_{k} (s   ; z_k  )  &=&   a^{i*}_{0_k}  (s)\mbox{ for all }  s \ge \tau_{k-1} \mbox{ when } z_k = (1, \tau_{k-1} ) , \\
 a^{i*}_{k} (s   ; z_k  )  &=&    0  \mbox{ for all }  s   \mbox{ when } z_k = (0, \tau_{k-1} ) 
.  \hspace{4mm}  \mbox{ \eop}
\end{eqnarray*} 
\end{cor}
By further using 
Theorem \ref{Thm_threshold_policy}, the (zero-starting) BR policies  $a^{i*}_{0_k} (\cdot )$ can be chosen to be threshold policies;  in other words
  any  best response strategy can be described completely using $M$-thresholds, say call them $\theta^{i*}_1, \cdots, \theta^{i*}_M$. 
  This implies the existence of an MT-strategy $(\theta^{i*}_1, \cdots, \theta^{i*}_M)$  among the BR strategies against any given strategy profile of opponents,
which completes the proof of  Theorem \ref{thm_monotonocity_of_thrshold}.  \eop}{}
\Cmnt{

 \subsubsection*{Last part of Theorem \ref{thm_monotonocity_of_thrshold}}


\textbf{Case 1:}\par
{\color{red}If $c_1=c_2 =c_3\dots=c_{M-1}=0$ and $c^i_M>\nu$}
Note that
$$P\big(\{\tau_{k}\leq \theta_k\}\cup \{\tau_{k}>\theta_k\}\big)=1 $$
\\
where the event $\{\tau_{k}>\theta_k\}$ means the agent never gets the lock $k$.
\\


If $\theta_{k+1}\leq \theta_k$ , then $\tau_k$ will be greater than $\theta_{k+1}$ with non-zero probability ($e^{-\beta \theta_{k+1}}$).
\\ Lets consider the case where $\tau_k>\theta_{k+1}$, 
even if it gets the $k-th$ lock, it will not get the $k+1^{th}$ lock, and hence will get zero reward, with additional cost $\nu \beta (\theta_k-\theta_{k+1}))>0$ of trying after time $\theta_{k+1}$.
But if the agent would have stopped at $\theta_{k+1}$ or earlier , it will get zero reward with zero additional cost which contradicts the fact that strategy with $\theta_{k+1}<\theta_k$ is optimal. So the optimal thresholds should be monotonically increasing i.e., $\theta_{k}\leq\theta_{k+1}$ .\par
\textbf {Case 2:}{\color{red} Under the assumption that $c^i_k>\nu \ \ \forall k$ We have, $$\hspace{20mm}  \theta_1\leq\theta_2\leq\theta_3\dots\leq\theta_M \hspace{20mm} $$  
If the agent contacts the lock at $\tau_k$ and lets consider the case where $\tau_k>\theta_{k+1}$, in this case the expected reward from $\tau_k $ on wards  will be zero. But if $\theta_{k+1}\ge\theta_k$  then the expected reward will be:
\begin{eqnarray*}
 &=&\hspace{-4mm}     \int_{\tau_{k}}^{\theta_{k+1}}   \Big (
   \big (c^i_k+v_{k+2}^i(\small (1, t\small );\pi^j) \big )-\nu\beta t\Big )\beta e^{-\beta t} dt \\
&& -\nu \beta (\theta_{k+1}-\tau_{k}) e^{-\beta(\theta_{k+1}-\tau_{k})}.\\
 &\ge&\hspace{-4mm}     \int_{\tau_{k}}^{\theta_{k+1}}   \Big (
  c^i_k-\nu\beta t\Big )\beta e^{-\beta t} dt  -\nu \beta (\theta_{k+1}-\tau_{k}) e^{-\beta(\theta_{k+1}-\tau_{k})}.
\end{eqnarray*}
It is easy to observe that the optimizer of the above expression is $\theta_{k+1}=T$, and  this holds for every $k$, which implies in the original expression ?

\textbf {Case 3:} if $c^i_k$ is a monotone increasing sequence, we have 
}
\eop

}
\vspace{-1mm}

\section{\textbf{Reduced Game $G$}}
\label{sec_ReducedGame} 

By  Theorem \ref{Thm_threshold_policy} any best response includes a threshold  or T-strategy   and further by Theorem \ref{thm_monotonocity_of_thrshold} at least one of the best response strategies is an MT-strategy.  By virtue of these results, one can find a NE (if it exists) in a much reduced game; the space of strategies in the  original game is infinite dimensional while that in the reduced game would be ${\cal R}^M$. We would show that there indeed exists a unique NE in the reduced game and analyze it by further reducing the dimension of the game to one.

\Cmnt{
 gives a very interesting result;  the sufficiency of threshold strategies. Because if the Nash equilibrium exists, the Nash equilibrium strategy profile will be made up of threshold strategies (it's possible that other strategy profile may also give the Nash equilibrium, i.e., we don't  have uniqueness yet). But to find Nash equilbrium, it is sufficient to work with threshold strategies.\par
 Theorem \ref{Thm_threshold_policy} and corollary \ref{cor_finite_collection} imply  that one threshold is sufficient to describe the best response policy corresponding to lock $k$, irrespective of $\tau_{k-1}.$
Thus a best response strategy can be described completely using $M$-thresholds,  $\theta_1, \cdots, \theta_M$, with $\theta_k$ being the maximum time till which  one continues to acquire $k$-th lock (in case it successfully acquired the first $(k-1)$-locks). The player need to choose only the best response ``time thresholds" for the  $M$ policies against the time thresholds of the opponent (due to symmetric nature of the game) rather than deciding the whole wave form (rate function) for the interval $[0,T]$. \par
}
We can reduce the problem to the  following game,
$G=\langle N,S, \Phi\rangle, $
where $N$ is the set of players as before,  $S$ the  set of strategies of each player is simplified 
%
to a  bounded set of  $M$ dimensional vectors (basically the set of MT-strategies):
$$
     S^i =  \big{\{}\underline{\theta}^i=(\theta^i_1,\theta^i_2,\dots \theta^i_M) ;  \  \theta^i_k \in [0,T] \  \forall k  
     \big{\}},
$$
 and the utilities $ \Phi = \{\phi^i\}$ are  now redefined by the following:
 \begin{eqnarray*}
 \phi^i ({\underline  \theta}^i ; {\underline  \theta}^j)  &=&
\phi^i ( {\underline  \theta}^i  ;  \ \theta^j_1)  =  \sum_{k=1}^M  E[ r_k^i( z_k^i, \theta_k^i;  \  \theta_1^j  ) ] \mbox{ where }  \\
{\underline  \theta}^j  & :=&   \{ { \underline \theta}^m \}_{m \ne i}   \mbox{ and } {  \theta}_1^j   :=  \{ {  \theta}_1^m \}_{m \ne i} ;
\end{eqnarray*}
 the above   objective function depends only upon the first thresholds of the  opponents ($\theta_1^j$)   because:   a) once the player gets the first lock successfully, the opponents have no incentives to try for the further locks, as (and after) their first contact is unsuccessful;  b) the success of any agent for first lock depends upon the failure of other agents for the same lock and hence on 
${  \theta}_1^j  $;  c) 
the  redefined terms (e.g., $r_k^i$, $\phi^i$ etc) depend only upon the MT-strategies, the thresholds of which are {\it defined using zero-starting optimal  policies, $\{a_{0_k}^{i*}\}$ of Corollary \ref{cor_finite_collection}}; and d) thus the thresholds do not depend upon $\tau^i_1$, the contact time of the first lock. 

  The simplified expressions under these special strategies are provided below.  
Recall the $k$-th component of the vector $(i.e., \theta_k^i)$ represents time threshold till which  one should attempt to contact $k$-th lock  with full intensity $\beta^i$, if 
 the previous   contact is before $\theta_k^i$;
 otherwise  one would not  attempt  for the next lock.   Hence  the cost under MT strategies simplifies to (as in \eqref{cost} and by rewriting the last term as an appropriate integral in equality `a'):

 \vspace{-4mm}
{ \begin{eqnarray*} 
r_k^i(z_k^i, \theta_k^i;  \theta^j_1) \hspace{-2mm} &=&\hspace{-1mm} 0 \mbox{ if $\tau^i_{k-1} \ge \theta_k^i$, else it equals }\\
r_k^i(z_k^i, \theta_k^i;  \theta^j_1)\hspace{-3mm}  &=& \hspace{-3mm} -\nu \beta (\theta_k^i-\tau^i_{k-1})e^{ - \beta^i (\theta_k^i -\tau^i_{k-1}) }\\
&&\hspace{-13mm} + \int_{\tau^i_{k-1}}^{\theta_k^i} \hspace{-1mm} ( c^i_k \eta_k^i (s) - \nu \beta^i   (s-\tau^i_{k-1})  )\beta^i e^{ - \beta^i (s-\tau^i_{k-1}) }  ds \\
& \stackrel{a}{=}&   \hspace{-3mm}  \int_{\tau^i_{k-1}}^{\theta_k^i}( c^i_k \eta_k^i (s) -\nu  )\beta^i e^{ - \beta^i (s-\tau^i_{k-1}) }  ds \mbox{, with }  \\
\eta_k^i  (s )\hspace{-1mm} &=& \hspace{-2mm} \indc  {k = 1} e^{ -   \sum_{m \ne i} \beta^m  ( s \wedge \theta^m_1 )}  +  \indc { k  > 1}  .  \label{eqn_opt}
\end{eqnarray*}}
\Cmnt{
In the above $\eta_k^i = \eta_k^i(s-\tau_{k-1};\Gamma_{\theta^j})$
 All the players choose the vector $\underline{\theta}$ at the beginning of the game, i.e., at $t=0$. Such a strategy profile means the player will try to contact the first lock with full potential ($\beta^i$, for the ease of notation we denote it simply by $\beta$) till time $\theta_1$, if the player is unable to contact till $\theta_1$ then the player will stop trying, and be silent for rest of the game. If the player is contacts the lock in time $\theta_1$ but the contact is unsuccessful, then also the player will stop trying, because the player will have no incentives for trying further, but if the player contacts the lock before $\theta_1$ and the contact is successful, then the player continues his search and starts looking for lock 2 with full potential till time $\theta_2$. If the player gets  he first lock successfully, the opponents become silent as they have no incentive to try further. For any lock $k$ if it contacts the lock before $\theta_k$, {\color{red} and if $\theta_{k+1}>\tau_{k}$} it starts looking for next lock,  if he is unable to contact the k-th lock before $\theta_k$  {\color{red}or if $\theta_{k+1}\le\tau_{k}$} he stops and remains silent for rest of the game.
 }

Once again we start with BR analysis, and  BR (of agent $i$) will be the maximizer of the following  objective function 

\vspace{-4mm}
{\small\begin{eqnarray}
\Upsilon_1^{i*}(z_1; \theta^j_1) =\max_{{\underline  \theta}^i   } \phi^i({\underline  \theta}^i ;  \theta^j_1)  
= \max_{ \{\theta_1^i, \cdots, \theta_M^i\}  }  \sum_{k \ge 1} E[ r_k^i( z_k^i, \theta_k^i; \theta^j_1 )].  \nonumber   \end{eqnarray}}
By applying  Theorem \ref{thm_monotonocity_of_thrshold} for finding best response against  MT-strategies of the  opponents   ($\pi^j := \{ { \underline \theta}^m \}_{m \ne i} )$)  we have:
  \begin{eqnarray*} \sup_ { \pi^i =  \{ a_1 (z_1) \} \cdots,  \{ \a_M (z_M) \} } \phi^i  (a_k \cdots, a_k;  z_1; \pi^j  )  
  =  \max_{\theta'_1 \cdots, \theta'_M}   \phi^i({\underline  \theta}^i ; {\underline  \theta}^j) ,  \end{eqnarray*}
  because  the optimal strategies can be chosen to be MT-strategies.  
  Now we apply DP equations to obtain the following, which  is further  simplified (in the second equation) 
 by choosing  an MT-strategy as the optimal strategy  for $\Upsilon^{i*}_2$ (once again by Theorem~\ref{thm_monotonocity_of_thrshold}): 
  
  \vspace{-4mm}
{\small \begin{eqnarray}
\Upsilon_1^{i*}(z^i_1; \theta^j_1)
 &=&  \max_{\theta^i_1} \gamma^i  (\theta^i_1;  \theta^j_1) \mbox{ where, }  \label{Eqn_further_reduced_game_util}    \\
\gamma^i  (\theta^i_1;  \theta^j_1) &:= &\hspace{-2mm}   \hspace{-2mm}  \int_{0}^{\theta^i_1}  \hspace{-1mm}  \Big (  \big ( c^i_1+\Upsilon_2^{i*} (t_1)  \big ) \eta_1^i (t_1) - \nu   \Big )\beta^i e^{-\beta^i t_1}dt_1 
\nonumber \mbox{ with }   \\
 \Upsilon_2^{i*}(t)  \hspace{-3mm} &:=& \max_{ \{\theta_2^i, \cdots, \theta_M^i\}  }  \sum_{k \ge 2} E[ r_k^i(z_k^i, \theta_k^i )   |  z^i_2 = (1, t)] .  \hspace{-5mm}\label{Eqn_Upsilon_two} 
 \end{eqnarray}}
 We now proceed with analysing the above BR and then the reduced game, as a first step,  we obtain the structural properties of  $\Upsilon_2^{i*}(t) $ (proof is in Appendix R\TR{}{,  \cite{TR}}): 
\Cmnt{
 By Corollary \ref{cor_finite_collection}, the maximizers defining  $\Upsilon_2^{i*}$ exists for all $t$ and  coincide in the intersecting intervals  (as  in the statement of the Corollary) for different values $t$;
 in other words we will see that these maximizers do not depend upon the strategies of the opponents. 
 In all,  we can see that   $\Upsilon_2^{i*}$ can be seen as  the cost under a   given (unique and maximizing)  sequence of $M-1$ thresholds $\{\theta^{i*}_k\}_{k \ge 2}$, that it is strictly  decreasing in time for all $t \le \theta_2^{i*}$ and all these are obtained in the following Theorem}
\Cmnt{@@@@@@@@@\\
We have a unique optimizer for each $k$ as a function $t_{k-1}$ as given below, which are defined backward recursively:
\begin{eqnarray}
\theta_k^{i*} (t_{k-1} )  &=& 
\left \{ 
\begin{array}{cccc}
   t_{k-1}      &  \mbox{ if }       t_{k-1}    >   \theta_k^{i*}     \\
\theta_k^{i*} &  \mbox{ else,   where  } 
   \end{array} \right . \nonumber  \\
 \theta_k^{i*}  &:=&     \inf_{t \ge 0 }   \{   c^i_k+\Upsilon_{k+1}^{i*} (t)        \le   \nu  \}   \wedge T.
\end{eqnarray}
The cost to go starting from $t_{k-1}$ is strictly decreasing for all $t_{k-1}  < \theta^{i*}_k$, after which it remains at 0: 

\vspace{-4mm}
{\small\begin{eqnarray*}
\Upsilon_M^{i*}(t_{M-1}) &=& 1_{c^i_M \ge \nu} e^{ \beta t_{M-1} } \ \int_{ t_{M-1}}^{ T} (c^i_M - \nu )\beta e^{-\beta  s  }d s \\
&&
\mbox { and for }  k < M \\
\Upsilon_{k}^{i*} (t_{k-1})   &=& \int_{t_{k-1}}^{  \theta^{i*}_k} \left [  c^i_k+\Upsilon_{k+1}^{i*} (t_k )  
   - \nu      \right ]
\beta e^{-\beta t_k  }d t_k
\end{eqnarray*}}
@@@@@@@@@}
\begin{thm}
\label{Thm_Upsilon_two}
Define the following    backward recursively:

\vspace{-4mm}
{\small \begin{eqnarray}
 \theta_M^{i*}  &=&  T \indc{ c^i_M > \nu}  \mbox{ and }  \nonumber \\
{\bar \Upsilon}_{M}^{i*} (t)   & = &   (c^i_M-\nu) \left ( 1 -  e^{- \beta^i (T-t) }  \right ) \indc{ c^i_M > \nu} ,   \nonumber  \\
&&\hspace{-18mm} \mbox{\normalsize and for any } 2\le    k < M \mbox{ (with $\emptyset$ -  null set)}   \nonumber \\
 \theta_{k}^{i*}  &:=&     \inf  \{ t \ge 0 :  c^i_{k}+ {\bar \Upsilon}_{k+1}^{i*} (t)        \le   \nu  \}   ,  \  \inf \emptyset := T,   \label{Eqn_theta_k_start_gt_2}  \hspace{7mm}\\ 
{\bar \Upsilon}_k^{i*}(t) &=& \indc{t <   \theta^{i*}_k }   \int_{t}^{ \theta^{i*}_k   } (c^i_k+ {\bar \Upsilon}_{k+1}^{i*} (s)-\nu)\beta^i e^{-\beta^i (s-t) }ds  . \nonumber  
\end{eqnarray} 
Then, }i) For any $k$, 
the   function ${\bar \Upsilon}_k^{i*}$ is   strictly decreasing with $t$ for all $t < \theta_k^{i*}$, after which it remains at 0. 
Further the co-efficients in \eqref{Eqn_theta_k_start_gt_2} are uniquely defined.  
\\
 ii)  The function 
$
\Upsilon_{2}^{i*} (\cdot) $ defined in (\ref{Eqn_Upsilon_two})     equals  $
{\bar \Upsilon}_{2}^{i*} (\cdot) $,   the former  is  optimized by   unique optimizers   $\{\theta_k^{i*}  \}_{k \ge 2}$ (defined  in  \eqref{Eqn_theta_k_start_gt_2}) and    is  strictly decreasing/remains at 0 as in (i).  
\eop
 
\end{thm}
\medskip
\Cmnt{

\newpage
The function $\Upsilon_{2}^{i*} (\cdot) $ defined in (\ref{Eqn_Upsilon_two})  is an example of  $M-1$ stage markov decision process. This problem can be solved using DP equations, with details as follows:
define 

\Cmnt{
To prove the uniqueness of the threshold policy, we will prove that there is a unique time threshold which is optimal, and the uniqueness of time threshold will imply the uniqueness of threshold policy. Note that the optimal time threshold for any lock $k$ is independent of the contact epoch $\tau_{k-1}$ of the previous lock,  and is given by:
\begin{eqnarray*}
\theta^{i*}_k&\in&  \arg \max_{\theta_k} \int_{0}^{\theta_k} (c^i_k+\Upsilon_{k+1}^{i*} (t)-\nu)\beta e^{-\beta t}dt.\\

with $\Upsilon_{k+1}^* (t)  $  defined recusrively backwards as below:
\begin{eqnarray*}
\Upsilon_{M}^* (t)   & = &\max_{\theta_M} \int_{0}^{\theta_M-t} (c^i_M-\nu)\beta e^{-\beta s}dss\\
\Upsilon_k^{*}(t) &=&\max_{\theta_k} \int_{0}^{\theta_k-t} (c^i_k+\Upsilon_{k+1}^* (t+s)-\nu)\beta e^{-\beta s}ds
\end{eqnarray*}
}

Note that  $$\theta_k^{i*}\in\arg\max_{\theta_k} \int_{0}^{\theta_k} (c^i_k+\Upsilon_{k+1}^{i*} (t)-\nu)\beta e^{-\beta t}dt
 $$
To prove the uniqueness of  optimizer $\theta^{i*}_k$, we will prove the function 
$$ \int_{0}^{\theta_k} (c^i_k+\Upsilon_{k+1}^{i*} (s)-\nu)\beta e^{-\beta s}ds$$  
  is a strict concave function on the interval $\theta_k\in[0,T]$ for every $k\ge 2$. To prove the concavity, we need to prove the derivative of this function with respect to $\theta_k$,
   $$ (c^i_k+\Upsilon_{k+1}^{i*} (\theta_k)-\nu)\beta e^{-\beta \theta_k}$$
is strictly decreasing for $\theta_k\in[0,T]$. Hence we need to prove that the function $\Upsilon_{k+1}^{i*} (t)$ is a strictly decreasing  function of  $t\in[0,T]$ as $c^i_k,  \nu, \beta $ are constants  and $e^{-\beta t}$ is a strictly decreasing funtion of $t$.

We will prove this using induction, for $k=M$ we have, 	
\begin{eqnarray*}
\Upsilon_M^{*}(t) &=& \int_{ 0}^{T-t} (c^i_M - \nu )\beta e^{-\beta s }d s \\
\end{eqnarray*}
Observe that this is a strictly decreasing function of $t\in[0,T]$. Assume the function $\Upsilon_i^{*}(t)$ is strictly decreasing function till  is true for $i=M, M-1\dots k+1$ and now, if we prove this holds true for $i=k$ we are done. For $i=k$,
\begin{eqnarray*}
\Upsilon_k^{*}(t) &=&\max_{\theta_k} \int_{0}^{\theta_k-t} (c^i_k+\Upsilon_{k+1}^{i*} (t+s)-\nu)\beta e^{-\beta s}ds\end{eqnarray*}
 by the induction hypothesis, $\Upsilon_{k+1}^{i*} (t) $ is a strict decreasing function of $t$ which implies $\Upsilon_{k+1}^{i*} (t+s) $ is also a strict decreasing function of $t$. We have the integrand is a decreasing  function of $t$, and so is the integral, it implies $\Upsilon_{k}^{i*} (t)$ as a decreasing function of $t$. And hence, we  have  $\Upsilon_{k}^{i*} (t) $ for every k is a strict decreasing function of $t$ till it becomes zero.\par 
 This proves the derivative is strictly decreasing in the interval $[0,T]$ which implies the function $\Upsilon_{k} (t)$ is a strict concave function in $[0,T]$ and hence will have a unique maximizer, i.e., 
 $$  \theta_k^{i*} =    \inf_{t \ge 0 }   \{   c^i_k+\Upsilon_{k+1}^{i*} (t)        \le   \nu  \}   \wedge T.
$$
\Cmnt{
Now we will prove the uniqueness of the optimizers (optimal time thresholds) as follows: we have for any k the optimal threshold is given by  $\theta^{i*}_k$
\begin{eqnarray*}
\theta^{i*}_k&\in&  \arg \max_{\theta_k} \int_{0}^{\theta_k} (c^i_k+\Upsilon_{k+1}^{i*} (t)-\nu)\beta e^{-\beta t}dt.\\
\end{eqnarray*}
On maximizing the above equation with respect to $\theta_k$ we have;
$$\theta_k^{i*}  =    \inf_{t \ge 0 }   \{   c^i_k+\Upsilon_{k+1}^{i*} (t)        \le   \nu  \}   \wedge T.$$
where $\Upsilon_{k+1}^{i*} (t)$ is a strict decreasing fuction of $t$ till it becomes zero, which implies $\theta_k^{i*}$ is unique. So, for every $k\in \{2,3,\dots,M\}$ we have the unique optimizer. For k=1, we have 
\begin{eqnarray*}
\theta^{i*}_1
&\in &  \arg \max_ {\theta_k} \int_{0}^{\theta_1} \left [   ( c_1 + \Upsilon_2^{i*} (t)  ) P(\Psi > t)    
   - \nu      \right ]
\beta e^{-\beta t_1}dt  
\end{eqnarray*}
in the above equation the integrand is a strictly  decreasing function of $t$, which implies it is a strict concave function which implies the uniqueness of $\theta^{i*}_1$.
Note that the control for the $k-th$ lock starts at $\tau_{k-1}$ and if $\tau_{k-1}$ is greater than $\theta^{i*}_k$, it means the agent will not try for the $k-th$ lock and will stop at $k-1$th lock. Hence, the optimal threshold can be written as a function of $\tau_{k-1}$ as follows:
\begin{eqnarray*}
\theta_k^{i*} (t_{k-1} )  &=& 
\left \{ 
\begin{array}{cccc}
   t_{k-1}      &  \mbox{ if }       t_{k-1}    >   \theta_k^{i*}     \\
\theta_k^{i*} &  \mbox{ else,   where  } 
   \end{array} \right . \nonumber  \\
\end{eqnarray*}}\eop

 \Cmnt{

 \subsubsection{Utility}
The expected utility of player $i$ when it chooses  $\underline\theta =(\theta_1,\theta_2,\dots \theta_M) $ and the opponents choose  $\underline\psi =(\psi_1,\psi_2,\dots \psi_M) $ 

\begin{eqnarray*}
\label{Eqn_Objectivei}
\phi^i (\theta^i, \theta^j ) =  \sum_{k=1}^M  E[ r_k^i(\ z_k^i, \theta_k^i;  \  \theta^j  ) ].
\end{eqnarray*}
For $M$ lock problem, if $\beta$ is the maximum rate of contact of the tagged player and $\beta_j=\sum_{m \ne i}\beta_m$  is the maximum possible rate of contact of the opponents, then we find the best response against threshold strategy of the opponents as follows:\par }
}

  In the view of the above discussions, the game breaks  into two problems. a) an optimization problem to find the optimal thresholds from the second lock onwards, which  is analysed in Theorem \ref{Thm_Upsilon_two};    b) a further reduced  one dimensional game with utilities given by  $\{\gamma^i\}$ of  \eqref{Eqn_further_reduced_game_util},  where each player (say player $i$) has to choose threshold ($\theta^i_1$),  used  for searching  the first-lock keeping in view of 
$\Upsilon_2^{i*}$   and the  strategies of the opponents.

\Cmnt {

\par One needs to solve these two problems to solve the game, and as we know by Corollary \ref{cor_finite_collection} that  the optimal time threshold for any $k$ can be given independent of $\tau_{k-1}$, we first solve the optimization problem, and obtain the optimal time thresholds $\theta^{i*}_2 \cdots, \theta^{i*}_M$ (for second lock on wards) for player $i$ by considering that the agent $i$ is successful with first lock at time $t = 0$ (because these parameters will remain unchaged for any $\tau_1$). We can find these optimal thresholds by backward induction in the following manner:
\begin{eqnarray*}
\theta^{i*}_M&=& \arg \max_{\theta_M} \int_{0}^{\theta_M} (c^i_M - \nu \beta t)\beta e^{-\beta t}dt- \nu \beta \theta_M e^{-\beta \theta_M }\\
&=& \arg \max_{\theta_M} \int_{0}^{\theta_M} (c^i_M-\nu)\beta e^{-\beta t}dt\\
\end{eqnarray*}
\begin{eqnarray*}
\Upsilon_M^{*}(t_{M-1}) &=& \max_{\theta_M}\int_{0}^{\theta_M- t_{M-1}} (c^i_M - \nu )\beta e^{-\beta t_M}dt_M \\
&=&\max_{\theta_M}\int_{ t_{M-1}}^{\theta_M} (c^i_M - \nu )\beta e^{-\beta (s-t_{M-1}) }d s \\
&=& e^{ \beta t_{M-1} } \max_{\theta_M}\int_{ t_{M-1}}^{\theta_M} (c^i_M - \nu )\beta e^{-\beta  s  }d s \\
\end{eqnarray*}
Thus $\theta_M^{i*} = T 1_{c^i_M \ge \nu} $ and thus   
\begin{eqnarray*}
\Upsilon_M^{*}(t_{M-1}) &=& 1_{c^i_M \ge \nu} e^{ \beta t_{M-1} } \ \int_{ t_{M-1}}^{ T} (c^i_M - \nu )\beta e^{-\beta  s  }d s \\
\end{eqnarray*}
Consider the case when $c^i_M > \nu$. 
In general for any $k$, by induction assume  $\Upsilon_{k+1}^{i*} (t)  $  is strictly  decreasing  in $t$.  Then for $k$-th lock we have 
\begin{eqnarray*}
\Upsilon_{k}^{i*} (t_{k-1}) = \max_{\theta_k} \left \{ 
\int_{0}^{\theta_k - t_{k-1} } \left [  c^i_k+\Upsilon_{k+1}^{i*} (t_k+t_{k-1})    
   - \nu      \right ]
\beta e^{-\beta t_k}dt_k \right    \},
\end{eqnarray*}by induction hypothesis the integrand is strictly decreasing, and hence is a strict concave function.  Thus we have unique optimizer, further clearly  the optimizer (if in interior) does not depend upon $t_{k-1}$ and is given by the following:
\begin{eqnarray}
\theta_k^{i*} (t_{k-1} )  &=& 
\left \{ 
\begin{array}{cccc}
   t_{k-1}      &  \mbox{ if }   c^i_k+\Upsilon_{k+1}^{i*} (t_{k-1})     <  \nu      \\
\theta_k^{i*} &  \mbox{ else,   where  } 
   \end{array} \right . \nonumber  \\
 \theta_k^{i*}  &:=&     \inf_{t \ge 0 }   \{   c^i_k+\Upsilon_{k+1}^{i*} (t)        \le   \nu  \}   \wedge T.
\end{eqnarray}
If in interior then 
$$
\Upsilon_{k}^{i*} (t_{k-1})   = \int_{0}^{\theta_k^{i*} - t_{k-1} } \left [  c^i_k+\Upsilon_{k+1}^{i*} (t_k+ t_{k-1})    
   - \nu      \right ]
\beta e^{-\beta t_k}dt_k,
$$which is strictly decreasing in $t_{k-1}$ as  $c^i_k+\Upsilon_{k+1}^{i*} (t_k + t_{k-1})    
   - \nu  > 0$  for all $t_k < \theta_k^{i*} - t_{k-1}  $ and because 
   $$
   c^i_k+\Upsilon_{k+1}^{i*} (t_k + t_{k-1})    
   - \nu   \mbox{ is also decreasing in } t_{k-1}.
   $$
The same logic is true even when $\theta_k^{i*} = T$. Thus by induction one can have unique maximizers for all $k$ and  $\Upsilon_{k}^{i*}$  is strictly decreasing function 
if  recursively the following is true:
\begin{eqnarray}
c^i_k+\Upsilon_{k}^{i*} (t_{k-1})     >   \nu  \mbox{ for all }  k,  t_{k-1}  \le \theta^{i*}_k  \mbox{ and }  
\end{eqnarray}

In otherwords we require the following condition which needs to be checked in backward inductively:
\begin{eqnarray}
c^i_k+\Upsilon_{k}^{i*} (\theta_k^{i*})     >   \nu  \mbox{ for all }  k  .
\end{eqnarray}
Or in other words this is true when 
$$
c^i_k > \nu \mbox{ for all }  k.
$$
In fact in this case $\theta_k^{i*} = T$ for all $k \ge 2.$

\newpage
\begin{eqnarray*}
\Upsilon_M^{*}(t_{M-1})=\max_{\theta_M}\int_{0}^{\theta_M-t_{M-1}} (c^i_M - \nu )\beta e^{-\beta t_M}dt_M
\end{eqnarray*}
and then,
\begin{eqnarray*}
\Upsilon_{k}^{i*} (t_{k-1}) = \max_{\theta_k} \left \{ 
\int_{t_{k-1}}^{\theta_k} \left [  c^i_k+\Upsilon_{k+1}^{i*} (t_k)    
   - \nu      \right ]
\beta e^{-\beta t_k}dt_k \right    \}
\end{eqnarray*}
with
\begin{eqnarray*}
\theta^{i*}_k&=&  \arg \max_{\theta_k} \int_{0}^{\theta_k} (c^i_k+\Upsilon_{k+1}^{i*} (t_k)-\nu)\beta e^{-\beta t}dt.\\
\end{eqnarray*}
After obtaining the optimal thresholds, we can find player $i$'s best response to find $\theta_1$ against any $\theta^j$ (which actually depends only upon $\theta^j_1 := \{\theta^{-i}_1 \}$) of other players
}

It is easy to observe that in the further reduced game,  the utilities  given in  equation (\ref{Eqn_further_reduced_game_util})    depend upon one dimensional strategy $\theta_1^i$ and one dimensional strategies of the opponents $\theta^j_1 $.
For any  $\theta^j_1 $  fixed, the partial  derivative of $\gamma^i$ with respect to $\theta_1^i$ is given by:
\begin{eqnarray}
\label{Eqn_derivative}
 \bigg (  \big ( c^i_1+\Upsilon_2^{i*} (\theta_1^i )  \big ) \eta_1^i (\theta_1^i ) - \nu   \bigg )\beta^i e^{-\beta^i \theta_1^i},
 \end{eqnarray}which is strictly decreasing by Theorem  \ref{Thm_Upsilon_two}, by the definition of $\eta_1^i$  and because $e^{-\beta \theta_1^i}$ is strictly decreasing. Thus the utility function is strictly concave in $ \theta_1^i$. Observe that the utility function is also continuous in $\theta^j_1$, therefore, we have a $n$-player concave game. 
Further by strict monotonicity of the derivative \eqref{Eqn_derivative}, the reduced game  satisfies {\it strict diagonal concavity}  given by \cite[equation (3.10)]{Rosen}.   Thus by \cite[Theorem 2]{Rosen}, we have unique NE for the reduced game.  It is easy to verify further details of the  following Theorem:
\begin{thm}\label{thm_unique_NE}
The unique NE is given by the sequence of thresholds (for second lock onwards) as given in Theorem \ref{Thm_Upsilon_two} (one sequence for each player)  and the first lock thresholds that simultaneously satisfy the following (for all $1\le i\le n$):

\vspace{-4mm}
{\small
\begin{eqnarray}
\theta_1^{i*}  =  \inf  \left \{ t : ( \Upsilon_2^{i*}  (t) + c^i_1 ) e^{ - \sum_{m \ne i} \beta^m  ( t \wedge \theta^{m*}_1 )}  
   \le  \nu   \right \} \wedge T . \label{Eqn_unique_NE} \hspace{3mm} \mbox{\eop}
\end{eqnarray}
 }
\end{thm} 

\subsection*{Some examples}

\TR{
\subsubsection*{Symmetric case  with large $T$}
When $M=2$   and consider the case with large $T$.  All symmetric agents.  By symmetry and uniqueness,  $\theta_k^{i*} = \theta_k^{*}$ for all $i \le n.$     
As $T \to \infty$,  we have that
$$
\Upsilon_2^* (t) = (c_2 - \nu)  (1 - e^{ - \beta (T -t) } ) \indc{ c_2 > \nu} \approx  (c_2 - \nu)^+  \mbox{ for }  t \ll T.  
$$
All the functions defining inifimum in \eqref{Eqn_unique_NE} are continuous and hence infimum is achieved and hence, 
\begin{eqnarray*}
\theta_2^{*}   & =&  T \indc{c_2 > \nu}   \mbox{ and  }  \\
\theta_1^*  &\approx & \max \left  \{ 0,      -\frac{1}{(n-1) \beta} \log \left ( \frac{ \nu } {  (c_2 - \nu)^+   + c_1  }  \right )  \right \}
\end{eqnarray*}

In fact when you substitute the approximate $\Upsilon_M^* (t) $ in  $\Upsilon_{M-1}^* (t)$, we obtain  (again with approximation as $\theta_M^* = T \indc{ c_M > \nu } $):
 
 \vspace{-4mm}
 {\small \begin{eqnarray*}
  \theta_{M-1}^* &=& T   \indc{(c_M-\nu)^+  + c_{M-1}- \nu > 0} ,  \\
\Upsilon_{M-1}^* (t) &=&  ( (c_M-\nu)^+  + c_{M-1}- \nu)  (1 - e^{ - \beta (T -t) } ) \\
& \approx  & ( (c_M-\nu)^+  + c_{M-1}- \nu)^+   \mbox{ for }  t \ll T .
 \end{eqnarray*}}
Progressing this way   for all $k >1$, define ${\bar c}_l^k := \sum_{l' = l}^k c_{l'}$.
 
 \vspace{-4mm}
 {\small \begin{eqnarray*}
\theta_k^* &=& T \indc {{\bar c}_k^{M^o_{k+1}}     \ge  (M^o_{k+1} - k + 1)   \nu} \mbox{ where }  M^o_{M}  = M -  \indc {c_M < \nu },  \\
M^o_k  &:=&      M^o_{k+1}   \indc { {\bar c}_k^{M^o_{k+1} }   \ge  (M^o_{k+1} - k + 1) \nu  } + (k-1)    \indc { {\bar c}_k^{M^o_{k+1}}   <  (M^o_{k+1} - k + 1) \nu  }   \\ 
&& 
\mbox{\normalsize and }  \\
\Upsilon_{k}^* (t)   &\approx &   \left (  {\bar c}_k^{M^o_{k+1} }   -   (M^o_{k+1} - k + 1)  \nu  \right )^+   \mbox{ for }  t \ll T,
\end{eqnarray*}}
and then 

\vspace{-4mm}
{\small \begin{eqnarray*}
\theta_1^*  &\approx & \max \left  \{ 0,      -\frac{1}{(n-1) \beta} \log \left ( \frac{ \nu } { { \bar c}_{1}^{M_2^o}    - (M_2^o-1) \nu}  \right )  \right \}   .
\end{eqnarray*}}

\subsubsection*{Monotone case  with large $T$} }{}

We consider an  asymmetric case,   in which the costs are monotone, i.e., without loss of generality assume\TR{}{\footnote{One can easily generalize a lot more (see \cite{TR}).}} $c^i_k \ge c^{i+1}_k $  for each    $i$, and that $c_M^n \ge M \nu$. 
  Also assume $\beta^i = \beta$ for all $i$. 
 We would derive the results that would be accurate for large $T$ and verify the same using numerical results. 
 
 With significantly large $T$, it is easy to observe (for all $i$):
 
 \vspace{-4mm}
 {\small
 $$
 \Upsilon_M^{i*} = (c_M^i - \nu)  (1 - e^{-\beta (T-t) } )  \approx   (c_M^i - \nu) \  \    \forall  t \ll T \mbox{ and }   \theta_M^{i*} = T .
 $$}
By substituting this   approximation in  {\small$\Upsilon_{M-1}^* (t)$ ($ \theta_{M-1}^{i*} = T$)}:
 
  \vspace{-4mm}
 {\small \begin{eqnarray*}
\Upsilon_{M-1}^{i*} (t) \approx  \hspace{-2mm}   \sum_{k=M-1}^M (c^i_k -  \nu)   (1 - e^{ - \beta (T -t) } ) 
 \approx \hspace{-2mm}    \sum_{k=M-1}^M (c^i_k -  \nu)     \mbox{ for }  t \ll T.  
 \end{eqnarray*}}
Progressing similarly   for all $k >1$, (with  ${\bar c}_k^{i} := \sum_{l = k}^M c^i_{l}$):
 
 \vspace{-4mm}
 {\small \begin{eqnarray*}
\theta_k^* \approx T  ,  \mbox{ and } 
\Upsilon_{k}^{i*} (t)   \approx    \left (  {\bar c}_k^{i}   -   (M  - k + 1)  \nu  \right )   \mbox{ for }  t \ll T,
\end{eqnarray*}}
Since the rewards of the players are monotone, we conjecture the corresponding $\theta_1^{i*}$ defining the NE  \eqref{Eqn_unique_NE} are also decreasing with $i$. We will derive the solution
of fixed point equations given in \eqref{Eqn_unique_NE}, which eventually  verifies the above conjecture. 
Further since  the functions defining inifimum in \eqref{Eqn_unique_NE} are continuous  (and strictly monotone), the   infimum is achieved and hence, 

\vspace{-4mm}{
\small $$\theta^{n*}_1=- \frac{1}{(n-1)\beta}\log\Big(\frac{\nu}{  {\bar c}^n_1  -(M-1)\nu }\Big),   \mbox{ if } \theta_1^{n *} \le \theta_1^{i*} \mbox{ for all }  i. 
$$}
 Now consider the player $n-1$,   repeating the same logic and substituting $\theta^{n*}_1$ derived in previous step we have:

\vspace{-4mm}
{\small  \begin{eqnarray*}
\theta^{(n-1) *}_1 &=&- \frac{1}{(n-2)\beta}\log\Bigg(\frac{\nu}{\big(  {\bar c}^{n-1}_1  -(M-1)\nu\big) e^{-\beta \theta^{n*}_1} }\Bigg)  \nonumber \\
&=&-  \frac{1}{(n-2)\beta}  \log\Bigg(\frac{\nu}{\big(  {\bar c}^{n-1}_1  -(M-1)\nu\big) } \Bigg) 
 -   \frac{  \theta^{n*}_1  } {  n-2   }. \hspace{8mm} 
\end{eqnarray*}}
As ${\bar c}^{n-1}_1 > {\bar c}^{n}_1$, it is indeed true that $\theta^{(n-1) *}_1  \ge \theta^{n *}_1 $. Progressing in exactly similar way, and verifying at every step the required monotonicity of 
$\{\theta_1^{j *}\}$, we obtain:

\vspace{-4mm}
{\small \begin{eqnarray}
 \label{Eqn_NE_large_T}
\theta^{(n-i)*}_1  =  \frac{ - \log\left (\frac{\nu}{   {\bar c}^{n-i}_1  -(M-1)\nu  } \right ) -  \beta  \sum_{j=n-i+1}^{n}  \theta^{j*}_1 }  { (n-i-1)\beta  }  .
\end{eqnarray}}
{\bf Remarks:} {\it
This solution is exact when $M=1$, thus we solved the problem completely for this case (solution matches with that in \cite{Mayank} when $n=2$); here the agents attempt to acquire the lock/destination without knowing if it already taken by others.}  
 For other cases it is an approximation, the accuracy is verified in the next section.

\section{Numerical Examples}


%
%
\Cmnt{
Consider the symmetric case, i.e.,  $\beta^i=\beta$  and $c^i_k=c_k\mbox{ for all } i\in N$. The uniqueness of NE given by Theorem \ref{thm_unique_NE} implies for any $k$, $\theta_k^i=\theta_k^j=\theta_k^*\mbox{ for all } i,j\in N$.\par
 The NE for $1$ lock and $n$ player symmetric game using Theorem  \ref{thm_unique_NE} is given by time thresholds $(\theta^*\dots\theta^*)$  with $\theta^*$ defined as:
\begin{equation*}
\theta^* := \left \{ \begin{array}{lllll}
   -\frac{  \ln ( \frac{\nu}{c_1} ) } { (n-1) \beta }   & \mbox{ if }   e^{-\beta (n-1) T}  \le  \nu\\
   T   &\mbox{ else.}
\end{array}      \right .
\end{equation*}
which also proves the conjecture given in \cite{Mayank} with $c_1=1$.}

We consider some numerical examples with an aim to reinforce the theoretical results. We also demonstrate that the approximation \eqref{Eqn_NE_large_T} is good even for moderate $T$. 
The first example considers symmetric case with moderate  $T= 8$ and the results are in Figure \ref{Fig_1}.  The other details are:  $M=5$, $n=4$, $\beta=1$,  $c_1=1,\ c_2=3,\ c_3=3, \ c_4=3, \mbox{  and }\ c_5=3 .$ We plot the NE for varying values of $\nu$. Theoretical results imply $\theta_k^* =T$ (for all $k > 1$) and we observe the same,  and, thus we plot only $\theta_1^{*}$.  
We computed the NE  by solving the fixed point equation \eqref{Eqn_unique_NE}, using fixed point iterates, we also plot the theoretical approximation \eqref{Eqn_NE_large_T}, and the two curves are indistinguishable (see Figure  \ref{Fig_1}). 
\Cmnt{
 We validated the results for $M$ locks and $n$ players symmetric  game; the uniqueness of NE ( by Theorem \ref{thm_unique_NE}) implies for any $k$, $\theta_k^i=\theta_k^j=\theta_k^*\ \forall \ i,j\in N$.  For large $T$ , we considered We took different values of $\nu$ ranging from $0.3$ to $2$; theoretical results imply $\theta_2^*=\theta_3^*=\theta_4^*=T$. We calculated $\theta_1^*$ using simulations and compared it with the theoretical result (see figure \ref{Fig_1}).  

}
\begin{figure}[H]
\begin{center}

\vspace{-6mm}

\begin{minipage}{4cm}
\vspace{2mm}
\hspace{-10mm}
\includegraphics[scale=0.2]{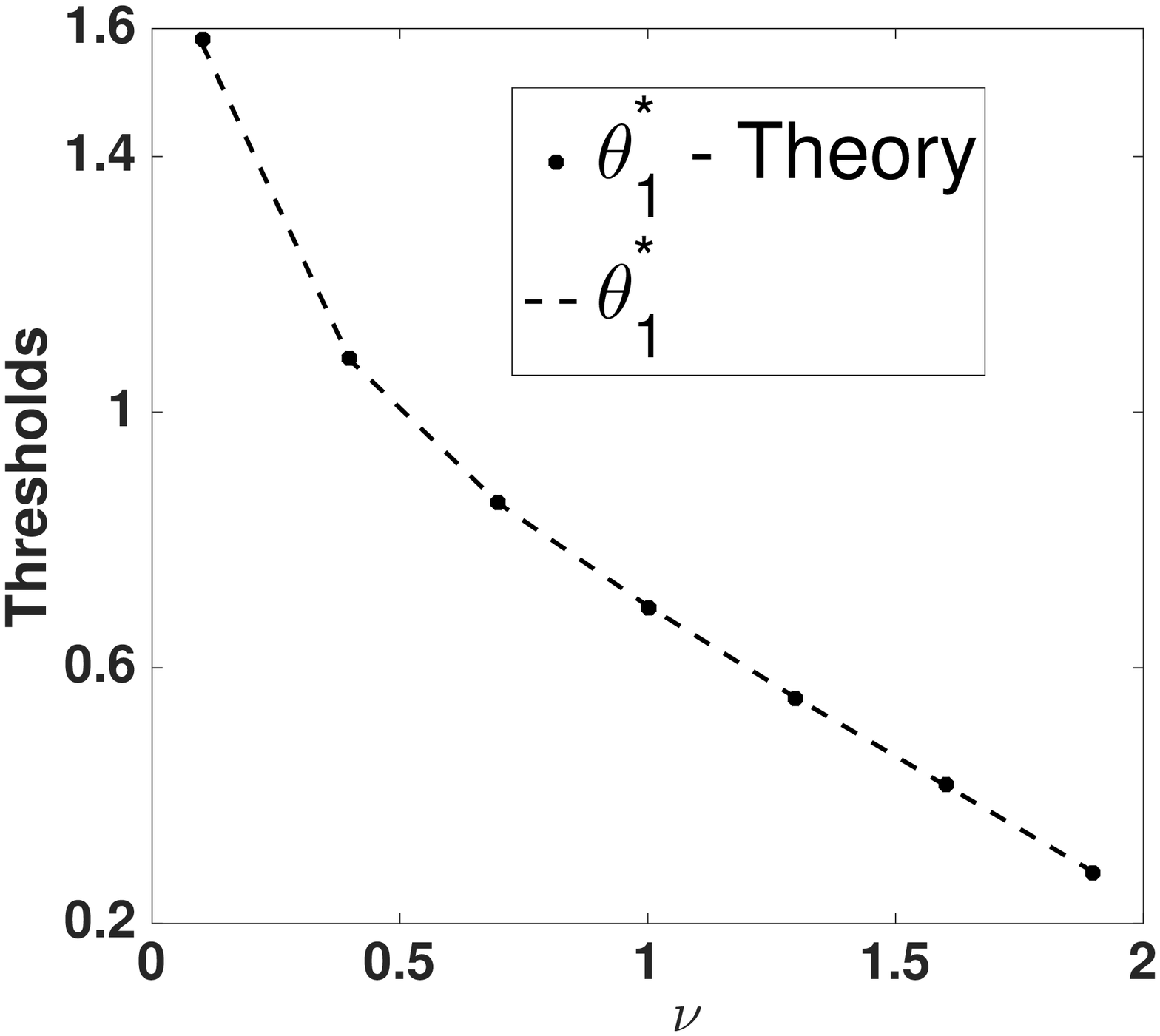}
\vspace{-13mm}
\caption{ Theoretical and numerical $\theta^*_1$ against $\nu$ (for large $T$) 
\label{Fig_1}}

\end{minipage}
\hspace{3mm}
\begin{minipage}{4cm}
\vspace{-8mm}
\includegraphics[scale=0.2]{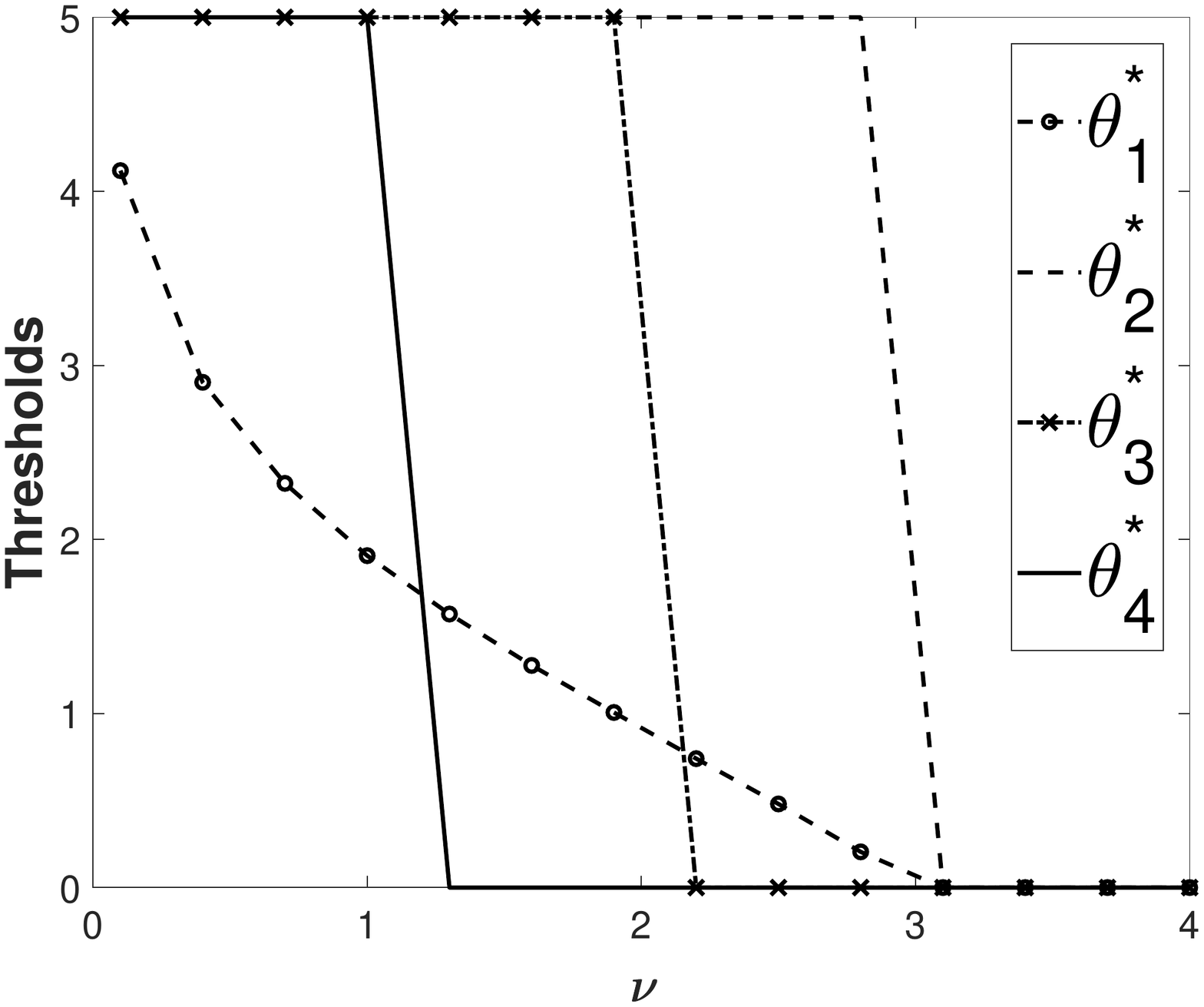}
\vspace{-17mm}
\caption{The NE thresholds against $\nu$ (for small $T$) 
\label{Fig_2}}

\end{minipage}
\end{center}
\vspace{-5mm}
\end{figure}
Further, to analyse the NE with small $T=5$, we consider a second case with  $M=4$ and  $n=2$ in Figure \ref{Fig_2};  we also set  $\beta=1$,   $c_1=4,\ c_2=3,\ c_3=2  \mbox{  and } c_4=1.$ 
For large values of $\nu$, the optimal thresholds are less than $T$ even for $k > 1$, thus even after acquiring the first (or a consecutive) lock the agent will not continue further if any of latter  locks are 
not acquired before the corresponding  thresholds.  \TR{
Also observe that, as the $\nu$ increases, the optimal threshold for $k=4$ becomes zero, while others are positive implying the agents will only attempt for three locks.  As $\nu$ increases further,  the optimal threshold for $k=3$  also becomes zero, while others are positive implying the agents will only attempt for two locks.
%

\begin{figure}[h]
\begin{center}

\vspace{-12mm}
\begin{minipage}{4cm}
\includegraphics[scale=0.2]{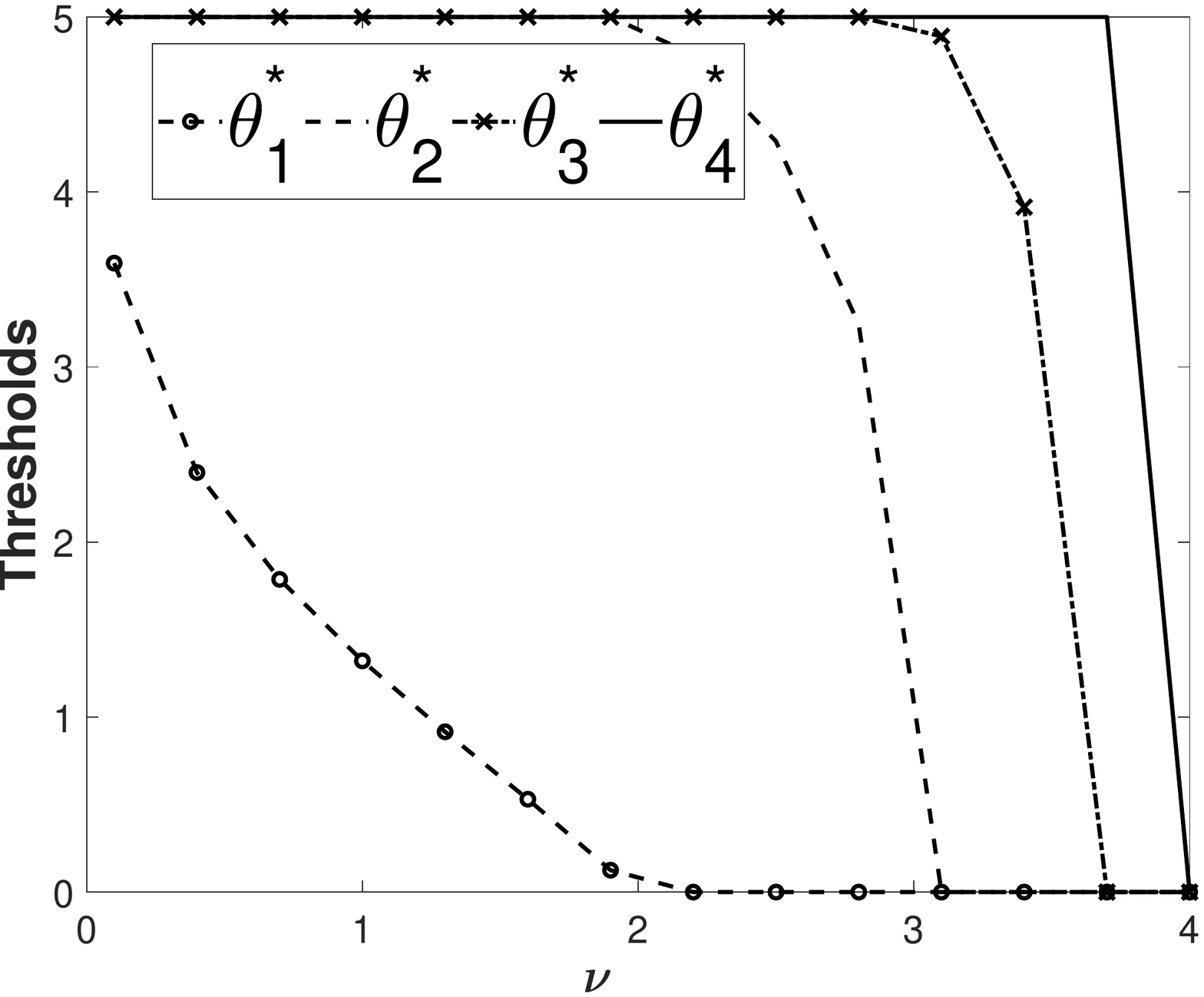}
\vspace{-13mm}
\caption{The NE thresholds for different vlaues of $\nu$ when $c_1= 1,  \  c_2= 2,\ \  c_3= 3,  \  c_4= 4   $
\label{Fig_3}}

\end{minipage}
\hspace{3mm}
\begin{minipage}{4cm}
\includegraphics[scale=0.2]{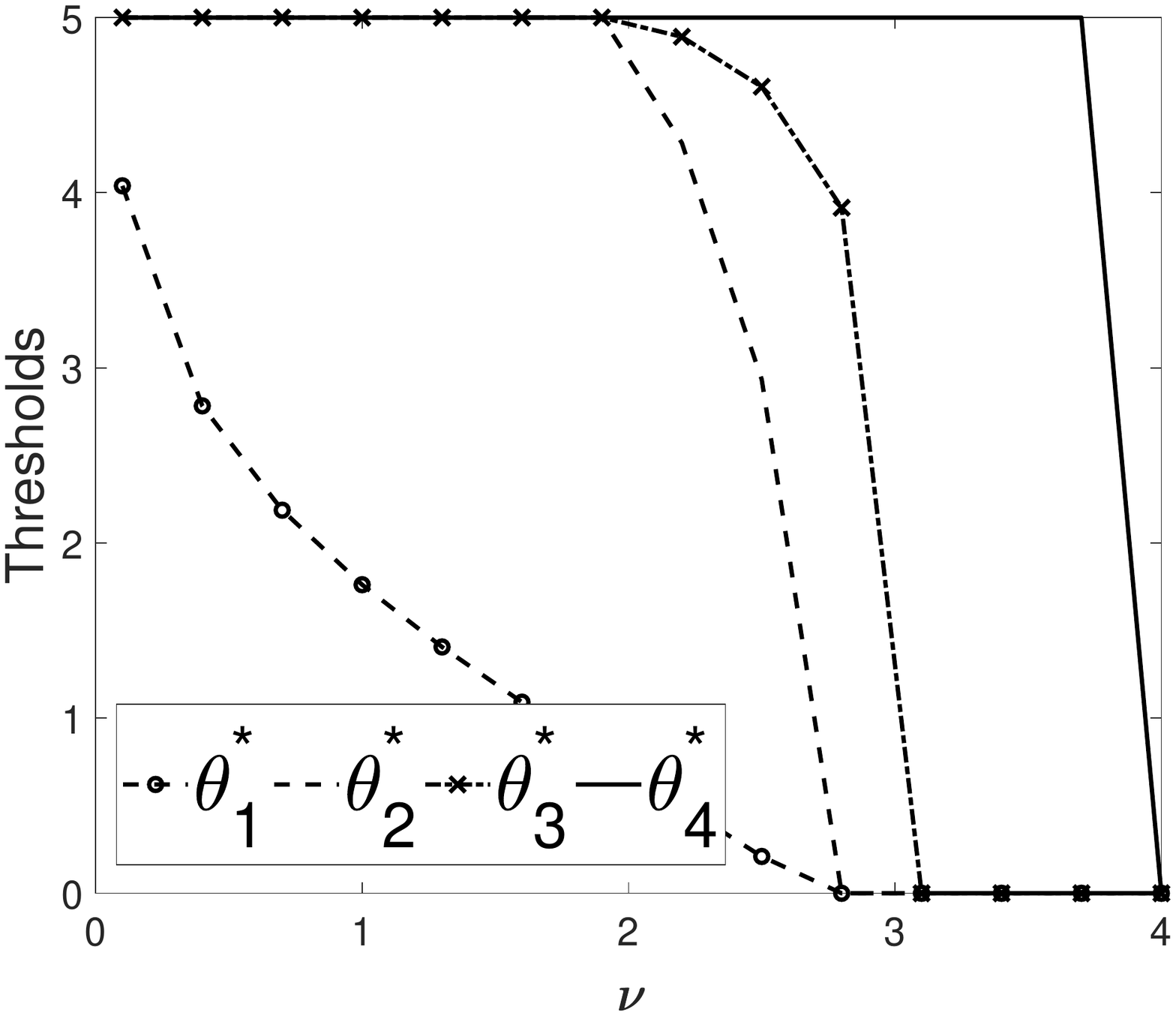}
\vspace{-13mm}
\caption{The NE thresholds for different vlaues of $\nu$ when $c_1= 4,  \  c_2= 2,\ \  c_3= 2,  \  c_4= 4  $
\label{Fig_4}}

\end{minipage}

\end{center}
\end{figure}

Further, to analyse the NE with small $T=5$, we consider more examples with  $M=4$ and  $n=2$ and  $\beta=1$. 
 In Figure \ref{Fig_3} we set  $c_1= 1,  \  c_2= 2,\ \  c_3= 3,  \  c_4= 4   $ and we see, when $\nu$ approaches 2.5, $\theta^*_1=0$, which implies the agents won't attempt the locks even though other thresholds are positive. In Figure \ref{Fig_4} we took the rewards associated with locks to be $c_1= 4,  \  c_2=2,\ \  c_3= 2,  \  c_4=4$. In all the examples, we observe that the optimal threshols are decreasing (non- increasing) with $\nu$.  
}{More examples will be explored in \cite{TR}.}

\Cmnt{
\begin{figure}[h]
\begin{center}

\vspace{-14mm}
\begin{minipage}{4cm}
\includegraphics[scale=0.2]{C_arr_4_4_4_4}
\vspace{-13mm}
\caption{The NE thresholds for different vlaues of $\nu$ when $c_1= 4,  \  c_2= 4,\ \  c_3= 4,  \  c_4= 4   $
\label{Fig_1}}

\end{minipage}
\hspace{3mm}
\begin{minipage}{4cm}
\includegraphics[scale=0.2]{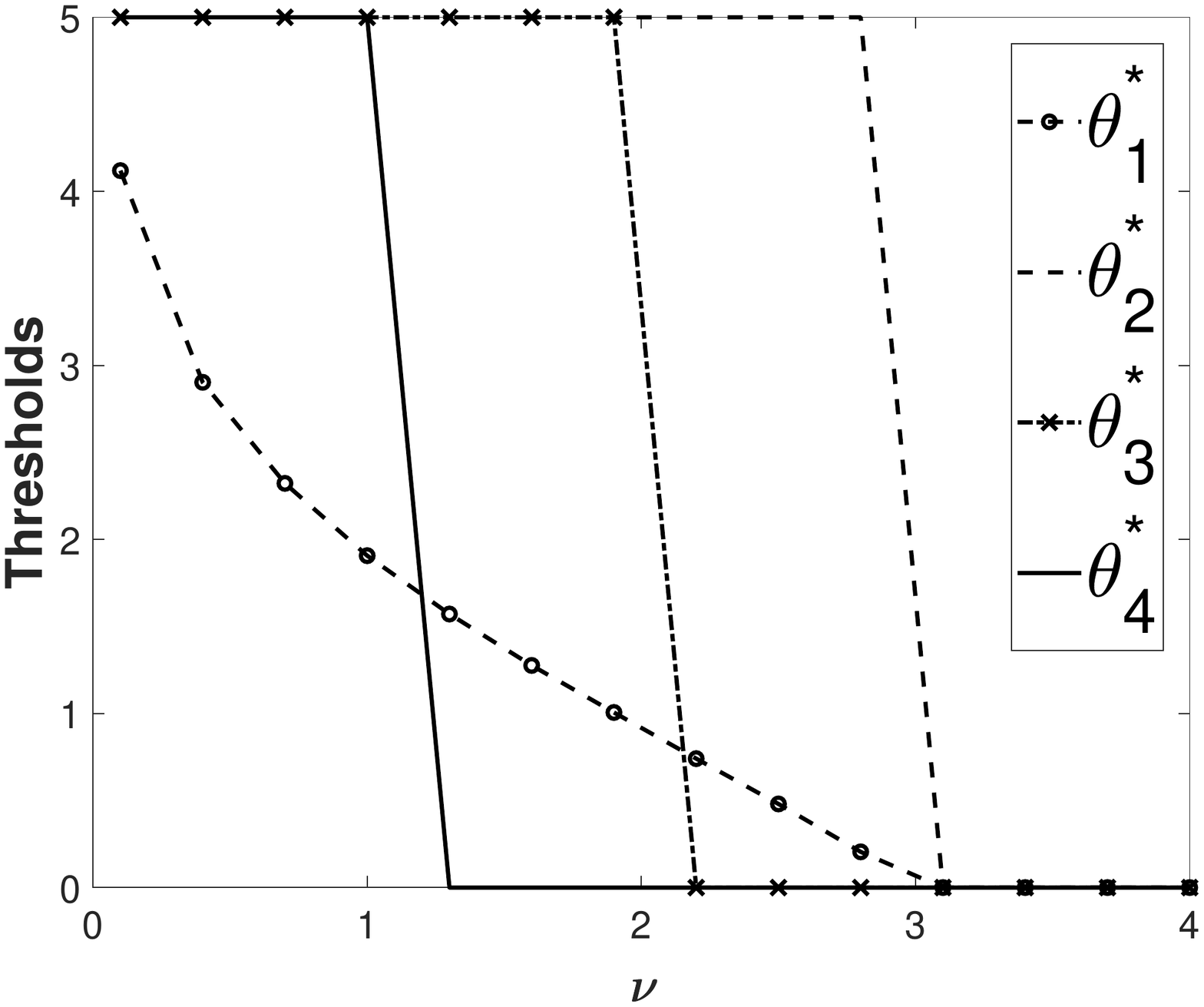}
\vspace{-13mm}
\caption{The NE thresholds for different vlaues of $\nu$ when $c_1= 4,  \  c_2= 3,\ \  c_3= 2,  \  c_4= 1   $
\label{Fig_4}}

\end{minipage}

\end{center}
\end{figure}
}

\Cmnt{
\subsection{
When others are silent}
In this case, $ P(\Psi > t) = 1$ for all $t$ and hence one needs to maximize 

{\small \begin{eqnarray*}
\Upsilon_{k}^{i*} (t_{k-1}) = \max_{\theta_k} \left \{ 
\int_{t_{k-1}}^{\theta_k} \left [  \Upsilon_{k+1}^{i*} (t_k)    
   - \nu      \right ]
\beta e^{-\beta t_k}dt_k \right    \}
\end{eqnarray*}}
the maximizer would be:
$$
\theta_k^{i*} (t_{k-1} ) = t_{k-1} \wedge \theta^{i*}_k
$$

 Since 
$\Upsilon^{i*}  $ is decreasing function, it is clear that
$$
\theta_k^{i*} = \sup_t \left   \{ \Upsilon_{k+1}^{i*} (t)    
   >  \nu \right  \}
$$
When $k = M$ (trying the last lock) then clear that $ \Upsilon_M^{i*} (t)   = c^i_M$ for all $t$ and hence 
$$
\theta_M^{i*} =  T 1_{ c^i_M > \nu }.
$$
And then 
$$
\phi^{i*}_M = (c^i_M - \nu )  (1 - e^{ - \beta^i T  } )   1_{ c^i_M > \nu }.
$$
When $k = (M-1)$, then (from second silence Theorem  or from above statement) we have
$$ \Upsilon_{M-1}^{i*} (t)   = c_{M-1} +  (c^i_M - \nu )  (1 - e^{ - \beta^i (T-t)  } )   1_{ c^i_M > \nu } \mbox{ for all } t$$ 
which is clearly decreasing with $t$.
Thus 
$$
(1 - e^{ - \beta^i (T-\theta_{M-1}^{i*} )  } )   =  \frac{  \nu -  c_{M-1}  } {c^i_M - \nu }
\mbox{ or }  $$
$$
\theta_{M-1}^{i*}  =  \frac{1}{\beta^i} \left ( T -  \log \left (  \frac{c^i_M - \nu } {  c^i_M - 2\nu +  c_{M-1}  } \right ) \right )
$$
We now have
\begin{eqnarray*}
 \Upsilon^{i*} (t)  & = & c_{M-2} + \sum_{l=0}^1  (c_{M -l}- \nu )  (1 - e^{ - \beta^i (\theta_{M-1}^{i*}  -t)  } )  \\ 
 &&   
-  (c^i_M - \nu )  e^{-\beta^i  T } ( \theta_{M-1}^{i*}  - t)
  \mbox{ for all } t \le \theta_{M-1}^{i*}  
\end{eqnarray*}
(Assuming some kind of monotonicity again)
Now, $\theta_{M-2}^{i*}$ satisfies

 {\small \begin{eqnarray*}
 \sum_{l=0}^1  (c_{M -l}- \nu )   e^{ - \beta^i (\theta_{M-1}^{i*}  - \theta_{M-2}^{i*} )  }   
 +  (c^i_M - \nu )  e^{-\beta^i  T } ( \theta_{M-1}^{i*}  - \theta_{M-2}^{i*} )  \\  =   \sum_{l=0}^2  (c_{M -l}-\nu)
\end{eqnarray*}}
And 

{\small
\begin{eqnarray*}
 \Upsilon^{i*} (t)  & = & c_{M-3} + \left (  \sum_{l=0}^2  (c_{M -l}- \nu )   -  (c^i_M - \nu )  e^{-\beta^i  T }  \theta_{M-1}^{i*} \right )  (1 - e^{ - \beta^i (\theta_{M-2}^{i*}  -t)  } )  \\ 
 &&   
-  \sum_{l=0}^1  (c_{M -l} -\nu )  e^{-\beta^i  \theta_{M-1}^{i*} } ( \theta_{M-2}^{i*}  - t)   \\
&& +   (c^i_M - \nu )  e^{-\beta^i  T }  \frac{ ( \theta_{M-2}^{i*}  - t)^2  }{2 }
  \mbox{ for all } t \le \theta_{M-2}^{i*}   
\end{eqnarray*}}

\newpage

\vspace{-4mm}
{\footnotesize{\begin{eqnarray}
  P_f &=& \int_{0}^{\theta_1}\Bigg(\int_{t_1}^{\theta_2}\Big(\int_{t_2}^{\theta_3}\dots\big(\int_{t_{M-1}}^{\theta_M}\beta e^{-\beta (t_M-t_{M-1})}dt_M\big)\nonumber\\
 && \hspace{10mm}\dots\beta e^{-\beta (t_3-t_2)}dt_3\Big)\beta e^{-\beta (t_2-t_1)}dt_2\Bigg)\nonumber\\
 &&\hspace{22mm} e^{-\beta_{j}(t_1\wedge\psi_1)}\beta e^{-\beta t_1}dt_1\nonumber
\end{eqnarray}}

 {\begin{eqnarray}
  P_f &=& \beta^M \int_{0}^{\theta_1} \int_{t_1}^{\theta_2} \int_{t_2}^{\theta_3}\dots \int_{t_{M-1}}^{\theta_M}  e^{-\beta t_M}dt_M  dt_{M-1}  \cdots dt_2 \nonumber\\ 
 &&\hspace{32mm} e^{-\beta_{j}(t_1\wedge\psi_1)}  dt_1\nonumber
\end{eqnarray}} }

{\color{blue}{\small
\begin{eqnarray*}
&=&\beta\int_0^{\theta_1} \Big(\beta^{M-2}
\int_{t_1}^{\theta_2} \int_{t_2}^{\theta_3}\dots \int_{t_{M-2}}^{\theta_{M-1}}\\
&& \hspace{20mm} \left ( e^{-\beta \theta_M }  - e^{-\beta t_{M-1}} \right )dt_{M-1}\\
&&\hspace{25mm}  \cdots dt_2\Big) e^{-\beta_{j}(t_1\wedge\psi_1)}  dt_1\\
&=&\beta\int_0^{\theta_1} \Big( \sum_{k=0}^{M-2} \mbox{\color{red}$($}(-1)^k \mbox{\color{red}$-\beta \nu \theta_{M-k-1})$} \beta^{M-k-2} e^{-\beta \theta_{M-k} }  \\
&& \hspace{10mm}
\int_{t_1}^{\theta_2}\dots \int_{t_{M-k-2} }^{\theta_{M-k-1}}   dt_{M-k-1} \cdots dt_2 \\
&&\hspace{15mm} + e^{-\beta t_1} (-1)^{M}\Big) e^{-\beta_{j}(t_1\wedge\psi_1)}  dt_1
\end{eqnarray*}
 } }
{\footnotesize
\begin{eqnarray*}
%
&=&\beta\int_0^{\theta_1} \Big( \sum_{k=0}^{M-2} (-1)^k \beta^{M-k-2} e^{-\beta \theta_{M-k} }  \\
&& \hspace{10mm}
\int_{t_1}^{\theta_2}\dots \int_{t_{M-k-2} }^{\theta_{M-k-1}}   dt_{M-k-1} \cdots dt_2 \\
&&\hspace{15mm} + e^{-\beta t_1} (-1)^{M}\Big) e^{-\beta_{j}(t_1\wedge\psi_1)}  dt_1
\end{eqnarray*}
  }
{\color{green}  
  All we need is the Lebegue measure/volume of the simplex 
  $$
 \left  \{  x_{i} \le \theta_i \forall 2 \le  i \le M-1, \ \ x_i \le x_{i+1} \forall 2\le  i \le M-2,   \ t_1 \le  x_2  \right \}  
  $$
  
  for M=4
\begin{eqnarray*}
  \{ t_1 \le x_2 \le \theta_2 , \  x_2 \le x_3 \le \theta_3  \} \hspace{-40mm} \\
 & =&  \{ t_1 \le x_2 \le \theta_2 , \  0  \le x_3 \le \theta_3  \}
  - \{ t_1 \le x_2 \le \theta_2 , \  0  \le x_3 \le x_2 \} \\
   & =& [t_1, \theta_2] \times [0, \theta_3] 
  -  1/2  [t_1,  \theta_2 ]^2   \\
\end{eqnarray*}

    $$
  \{ t_1 \le x_2 \le \theta_2 , \  x_2 \le x_3 \le \theta_3 \   \  x_3 \le x_4 \le \theta_4   \}
  $$
  }

and the expected cost will be


{\small{\begin{eqnarray}
 E[cost] &=& \beta^{M+1}\int_{0}^{\theta_1}\int_{t_1}^{\theta_2}\int_{t_2}^{\theta_3}\dots\int_{t_{M-1}}^{\theta_M} t_M e^{-\beta t_M} {dt_M}\nonumber\\
 &&\hspace{20mm}\dots dt_3{dt_2}e^{-\beta_{j}(t_1\wedge\psi_1)} dt_1\nonumber\\
 &&+ \beta^{M}\int_{0}^{\theta_1}\int_{t_1}^{\theta_2}\int_{t_2}^{\theta_3}\dots\int_{t_{M-2}}^{\theta_{M-1}} \theta_M e^{-\beta \theta_M} {dt_{M-1}}\nonumber\\
 &&\hspace{20mm}\dots dt_3{dt_2}e^{-\beta_{j}(t_1\wedge\psi_1)} dt_1\nonumber\\
&&+ \beta^{M-1}\int_{0}^{\theta_1}\int_{t_1}^{\theta_2}\int_{t_2}^{\theta_3}\dots\int_{t_{M-3}}^{\theta_{M-2}} \theta_{M-1} e^{-\beta \theta_{M-1}} \nonumber\\
 &&\hspace{20mm}{dt_{M-1}}\dots dt_3{dt_2}e^{-\beta_{j}(t_1\wedge\psi_1)} dt_1\nonumber\\
 &&\hspace{22mm}e^{-\beta_{j}(t_1\wedge\psi_1)}\beta dt_1\nonumber\\
 &&\hspace{15mm}\vdots\nonumber\\
  &&\hspace{15mm}\vdots\nonumber\\
  && + \beta^2 \int_{0}^{\theta_1}    \theta_2 e^{-\beta\theta_2}   
   e^{-\beta_{j}(t_1\wedge \psi_1)} dt_1\nonumber\\
  &&+\beta\theta_1 e^{-\beta \theta_1}\nonumber
\end{eqnarray}} }

{
\footnotesize
\begin{eqnarray*}
E[cost]&=& \sum_{k=1}^{M} \beta^{k} \theta_{k} e^{-\beta\theta_{k}}\int_{0}^{\theta_1}\int_{t_1}^{\theta_2}\\
&&\hspace{18mm}\dots\int_{t_{k-2}}^{\theta_{k-1}}dt_{k-1}\dots dt_2  e^{-\beta_{j}(t_1\wedge \psi_1)} dt_1\\
&&+ \beta^{M+1}\int_{0}^{\theta_1}\int_{t_1}^{\theta_2}\int_{t_2}^{\theta_3}\dots\int_{t_{M-1}}^{\theta_M} t_M e^{-\beta t_M} {dt_M}
\end{eqnarray*}}
{\color{blue}{Theorem statements}
\begin{itemize}
    
    \item A threshold policy  is always  in Best response  (Can we get uniqueness) 
    
    \item If uniqueness is done,  it is equivalent to work with threshold policies,  (try)
    
    \item Else, it is sufficient to work with threshold policies.. 
    
    \item  Monotonicity of thresholds 
    
    \item Define the Reduced game 
    
\end{itemize}
 
 We obtained the results of the reduced game through numerical simulations and discussed some interesting observations. 
}
}

\Cmnt{ 
\textbf{\textit{Resistance of the environment:}} \textit{The resistance of the environment for any lock $k$ is the probability that $kth$  contact will be unsuccessful (as function of time) by the tagged player depending on the state $z_k$ of the player. This is nothing but the probability that the other player would have already contacted the $kth$ lock at the time of contact of tagged player,i.e., 
 $$R_{z_k}(\cdot)=1-P^{z_k}_f(\cdot)$$
 where $P^{z_k}_f(\cdot)$ is the probability function which spans over $(\tau_k,T]$ and gives the probability of failure of other agent at any time $t$ in this interval.}\par
For any measurable function $f:[t, T] \to {\cal R}$ let $f_{[a, b]} : [a,b] \to {\cal R} $ represent a sub-function restricted to domain $[a,b] \subset [t, T]$, i.e., we have by notation:
$$
f(s) = f_{[a, b]} (s)  \mbox{ for all }  s  \in [a, b]. 
$$We use $f \equiv g$, when the two functions are the same and have the same domain.

Lets denote $a_{k,t}(\cdot)$ the control used for $k$th lock starting from time $t$ (it will span till the time $T$), also note $a_{k,t, [\tau,T]}$ denotes the part of the control $a_{k,t}(\cdot)$ from $\tau$ to $T$.

{\bf 
Recall that for ordered case, the resistance function is zero for stage $k >1$, i.e., 
$$
R_{z_k} \equiv 0, \mbox{ when }  z_k = (1, \tau) \mbox{ for some } \tau \le T. 
$$
The the agent would try only with $z_k = (1, \tau)$ (when $k > 1$) and hence the following lemma is only called either for $k = 1$ or for $k> 1$ with special states of the form, $z_k = (1, \tau)$}
}

{\Cmnt
{\color{red}
\begin{lemma}\label{lemma_bang_bang_policy}
Optimal (Best response) strategy of the player against any strategy of the opponent, is made up of bang bang policies, i.e., is a bang bang strategy.
$$\pi^*\in \Pi_\mathcal{B}.$$

\end{lemma}{}
Proof is in Appendix A.


\begin{lemma}\label{lemma_earlier_control}
For any lock $k$, consider an open loop policy $a(t)$ and construct a threshold policy $\Gamma_\psi$ such that 
$$\int_{0}^{T}a(t)dt=\int_{0}^{T}\Gamma_\psi(t)dt$$ then, we have the following
\begin{enumerate}
    \item The contact epoch $\tau_k$ using a threshold policy will be earlier with probability 1, i.e.,$$P(\tau_{\psi}\leq \tau_a)=1 \mbox{ and }P(\tau_{\psi}<\tau_a)>0.$$
    \item The probability of contacting a lock successfully (before the opponent contacts it) is more using threshold policy,  if $\eta$ denote the minimum of the contact epochs of the opponents, then
$$P\{\omega:\tau_{a}(\omega)\leq\eta(\omega)\} \le P\{\omega:\tau_{\psi}(\omega)\leq\eta(\omega)\}. $$
\end{enumerate}{}
\end{lemma}{}
Proof is in Appendix B. 
}}




\section*{Conclusions}
We considered acquisition games with partial, asymmetric information. Agents attempt to acquire $M$ destinations, the first one to contact a destination acquires it; destinations can be acquired only in a given order. When an agent  acquires a destination and if it has also acquired all the previous ones, it gets some reward. The agents are not aware of the acquisition status of others. It is possible that an agent continues its acquisition attempts,  while  the destination is already acquired by another agent. 
Thus we have a partial and asymmetric information game. 
 The agents control the rate of their Poisson clocks to contact the destinations;   they  also incur  a cost proportional to their rates of contact. We found NE of this asymmetric game by  reducing  it to a much simpler game such that the NE of the reduced game would also the NE of the original game; the original game has infinite number of state and has infinite dimensional actions.  We proved that a tuple of time-threshold policies form the unique NE of the reduced game. We also provided an algorithm  to compute this NE. We further have approximate closed form expressions for the NE.



\section*{Appendix   }


\noindent {\bf Proof of Theorem \ref{thm_existance}:}
The  DP equation for any $k$ can be rewritten~as:

\vspace{-4mm}
{\small 
\begin{eqnarray*}
u(s,x)\hspace{-3mm}& =& \hspace{-3mm} \sup_{ \ a (.) }   \left \{  \int_s^T  L_k(s',  x (s') ,   a(s') )  ds'   + g(x (T)  ) \right \} \mbox{  with }  \\
L_k(s',  x,  a)  \hspace{-3mm} &=& \hspace{-3mm}
  \left (  c^i_k + v_{k+1}^i( (1,s'); \pi^j) \big) \eta_k^i(s';\pi^j)     - \nu   x  \right ) a e^{-x}, \\
   \eta_{k}^{i}(s',\pi^j) \hspace{-3mm} & =  & \hspace{-3mm} e^{  - \sum_{m \ne i}    \bar{a}^m(s')} \indc{k= 1} + \indc{ k > 1}.
\end{eqnarray*}}
It is easy to solve \eqref{Eqn_DP} with $k=M$, because $v^i_{M+1} $ is defined to be 0 (details in  \cite{Mayank}, also discussed in Theorem \ref{Thm_Upsilon_two}): 
$$
v^i_{M} (1, s)    =    (c_M^i-\nu) \left ( 1 -  e^{- \beta^i (T-s) }  \right ) \indc{ c_M^i  > \nu} .
$$
   Observe that $v^i_{M} (1, s) $  is a differentiable function of $s$, hence clearly $L_{M-1}$ is Lipschitz continuous on $[t, T) \times [0, \beta T] \times [0, \beta]$. It is also bounded.  Further clearly the RHS of the  ODE and the terminal cost $g$ are all bounded and Lipschitz continuous. Thus by  \cite[Theorem 10.1 and the following Remark 10.1]{feller}  the value function $u(\cdot, \cdot)$ is unique viscosity solution which is Lipschitz continuous when $k= M-1$.  This implies 
$v^i_{M-1}(1, s)$ is     Lipschitz continuous in $s$, further it is also bounded. 
    This implies $L_{M-2}$ is also  Lipschitz continuous and bounded which proves the same for $v^i_{M-2}(1, s)$.  By backward induction on $k$, the part (i) is true. 
    
For part (ii) we apply the results of \cite{Roxin}. Towards this the optimal control problem can be converted into Mayer-type (finite horizon problem with only terminal cost) by usual technique of augmenting a new component to state which represents  \vspace{-4mm}
$$
y(s') := \int_{s}^{s'} L_k( {\tilde s},  x({\tilde s}),  a ({\tilde s}) )  d{\tilde s}
$$  and equivalently maximizing  $y(T) + g(x(T) )$. By part (i) all the required assumptions \cite[Assumptions (i) to (vii)]{Roxin} are satisfied, with compact control space  $U = [0,\beta]$, compact state space 
${\hat X} = [0, \beta T]$  (it is easy to verify that the state variable could be confined to this range);    assumptions (i)-(ii) are trivially satisfied; assumption (iii) follows by part (i); for assumption (iv)  one can actually bound by a constant independent of $(s', (x,y))$; easy to verify convexity  requirement of (vii), because for any given $(s', (x,y))$ the set in mention is an interval. 

  This proves part (iii).  \eop

\medskip

{\bf Proof of Lemma \ref{lemma_earlier_control}:} 
For any lock $k$, given any state $z^i_k$, consider the open loop policy $a_k^i(\cdot)$, we denote it as $a(\cdot)$ for the ease of notation. If $a(\cdot)$  is already a threshold type,  we will have nothing to prove. 
If not,  choose two  intervals  $[t_1,t_1+\delta_1]$ and  $[t_2,t_2+\delta_2]$   with $t_2 \ge t_1+\delta_1$  such that
$$\int_{t_1}^{t_1+\delta_1}a(t)dt<\beta \delta_1 \mbox{ and }  \int_{t_2}^{t_2+\delta_2}a(t)dt>0.$$
One can further ensure (by choosing appropriate end points) that
  $$\int_{t_1}^{t_1+\delta_1}a(t)dt + \int_{t_2}^{t_2+\delta_2}a(t)dt=\beta \delta_1. $$
  Now we
construct another policy, ${a'}(t)$ such that,
$$\int_{t_1}^{t_1+\delta_1}{a'}(t)dt=\beta \delta_1 \mbox{  and } \int_{t_2}^{t_2+\delta_2}{a'}(t)dt=0,$$ 
and on rest of the intervals the policy ${a'}(\cdot)$ matches completely with policy $a(\cdot)$. This new policy is basically constructed by shifting the mass from a later interval  $[t_2,t_2+\delta_2]$ to a former interval $[t_1,t_1+\delta_1]$ in policy $a(\cdot)$, note that if one can't find such intervals, it implies that the policy $a(\cdot)$ itself is  a Threshold policy and we will have nothing to prove. \par
Observe that  for all $t < t_1$ we have,
$$
 a(t) = {a'}(t) \mbox{ and hence  } \bar{a}(t) =\int_{\tau_{k-1}}^t a (s) ds = \bar{a}'(t) 
$$
similarly for all  $ t_1 <  t < t_2 $,
$\bar{a}(t)  < \bar{a}'(t) $ and   for all  $ t_2< t < t_2 + \delta
$  we have $  \bar{a}(t)  \le  \bar{a}'(t)$   and  for all  $  t  > t_2 + \delta$ we have   $  \bar{a}(t)  =   \bar{a}'(t)
$. \par
This implies, the time to contact the $k$-th lock with  policy $a(\cdot)$ denoted as  $\tau_a$ is stochastically dominated by that  with  policy $a'(\cdot)$ denoted as  $\tau_{a'}$\TR{ as explained below; consider the CDFs with both the policies; i) for any $x<t_1$,  
$$F_{a}(x) = Prob (\tau_a \le x) = 1-  e^{-\bar{a}(x) }    =   1-  e^{-\bar{a}'(x) }   =F_{a'}(x),$$ 
 ii) for any $x\in (t_1,t_1+\delta_1),$ we have $\bar{a}(t)<\bar{a}'(t)$ and so
 $$F_{a}(x)= 1-  e^{-\bar{a}(x) }    <   1-  e^{-\bar{a}'(x) }   =F_{a'}(x),$$
iii)  for any $x\in (t_2,t_2+\delta_2),$ we have $\bar{a}(t)\le\bar{a}'(t)$
 $$F_{a}(x)= 1-  e^{-\bar{a}(x) }    \le   1-  e^{-\bar{a}'(x) }   =F_{a'}(x),$$
and iv) for any $x\in [t_2+\delta_2,T],$ we have $\bar{a}(t)=\bar{a}'(t)$
 $$F_{a}(x)= 1-  e^{-\bar{a}(x) }    =   1-  e^{-\bar{a}'(x) }   =F_{a'}(x).$$
This proves the required stochastic dominance,
  $\tau_a \stackrel{d}{\le } \tau_{a'}$. }{ (the CDFs are compared, see  \cite{TR} for  details).}
%
%
\TR{Now, let $\tau^j=\min_{m\ne i} \tau^m $ where $\tau^m$ denotes the contact epoch of $m$-th agnet, and observe
$$\{\omega:\tau_{a}(\omega)\leq\tau^j(\omega)\} \subset\{\omega:\tau_{a'
}(\omega)\leq\tau^j(\omega)\}, $$
which is same as saying,
$$\{\mbox{success with policy } a \}\subset\{\mbox{success with policy } a' \}$$
 and hence, \vspace{-4mm}
 $$\hspace{10mm}  P_k^i(z^i_{k}; a; \pi^j)  \le  P_k^i(z^i_{k}; a'; \pi^j)  \mbox{ for any } k, z_{k} . $$}{Stochastic dominance implies (details in \cite{TR}) 
 $$\hspace{10mm}  P_k^i(z^i_{k}; a; \pi^j)  \le  P_k^i(z^i_{k}; a'; \pi^j)  \mbox{ for any } k, z_{k} . $$ }
Further observe that  the expected cost with policy $a(\cdot)$ (by change of variables)

\vspace{-4mm}
\TR{\begin{eqnarray*}
E[ {\bar a} (\tau_a ); \tau_a < T ] +  {\bar a} (T ) e^{-{\bar a}(T)} \hspace{-25mm} \\
& =& \int_0^T {\bar a} (t)  e^{- {\bar a} (t) } a(t) dt   +  {\bar a} (T ) e^{-{\bar a}(T)}  
\\
&=&   \int_0^{{\bar a}(T) }  x  e^{-x}  dx +  {\bar a} (T ) e^{-{\bar a}(T)}   =  1 -   e^{-{\bar a}(T)},  
\end{eqnarray*}which 
is the same as that using $a'$ because  ${\bar a}(T) = {\bar a}'(T)$.}

{\begin{eqnarray*}
E[ {\bar a} (\tau_a ); \tau_a < T ] +  {\bar a} (T ) e^{-{\bar a}(T)}& =&  1 -   e^{-{\bar a}(T)},  
\end{eqnarray*}which 
is the same as that using $a'$ because  ${\bar a}(T) = {\bar a}'(T)$. \TR{}{The detailed proof is provided in \cite{TR}.}}

%
\Cmnt{
For the part 2, we know the contact will be successful only if, the agent is able to contact the lock before the opponent. If the opponents are using the strategy $\pi^j$, and suppose using this strategy, it will be able to contact the lock at time $\tau^j$, note that $\tau^j$ is a exponentially distributed random variable with parameter $a^j=\sum_{m \ne i} { a}_k^m (t)$ .{\color{red} Note that there exists a $\delta>0$  such that $$P(\tau_a-\tau_{a'}>0)=\delta,$$let's consider this case. Note that $\tau^j\in (\tau_{a'}, \tau_a)$ with probability given below:
\begin{eqnarray*}
P(\tau_{a'}<\tau^j<\tau_{a})&=&F_{a^j}(\tau_a)-F_{a^j}(\tau_{a'})\\
&&= 1-  e^{-\bar{a}^j(\tau_a) }- 1+ e^{-\bar{a}^j(\tau_{a'}) } \\
&&=    e^{-\bar{a}^j(\tau_{a'}) }  -e^{-\bar{a}^j(\tau_a) } \\
&&>0
\end{eqnarray*}
 Note that the last inequality follows from the fact that the opponents are not silent in the time interval $(\tau_{a'}, \tau_a)$  }}
One can keep on improving the policy until it becomes a threshold policy. This completes  the proof. \eop

\par
{\bf Proof of Lemma \ref{Lemma_value_function_mono}:} 
 For any lock $k$, given any state $z^i_k$, the $k$-th stage DP equation can be re-written as the following optimal control problem (see equation (\ref{eqn_opt_ctrl})) 
 
 \vspace{-4mm}
 {\small\begin{eqnarray*}
u(t,0) =  \sup_{ \ a (\cdot) \in L^\infty}   \Big \{ \int_{s}^{T} \hspace{-3mm}  \big ( h_k^i (s')- \nu   x (s')\big )   
  a ( s') e^{-x(s')} ds'  +g(x(T )) \Big \}. 
\end{eqnarray*}}
%
  From the Dynamic programming principle \cite[Theorem 5.1]{feller}, we can rewrite it as follows;
 
 \vspace{-5mm}
 {\small \begin{eqnarray*}
u(t,0 )  &=& \hspace{-4mm}\sup_{ a_k^i \in L^\infty [t, \tau] }\Bigg\{ \int_{t}^{\tau} \big( h_k^i (s)- \nu   x (s)  \big )  a_k^i ( s)  e^{-x(s)}ds   \\
&& \hspace{18mm}+ u(\tau,x(\tau)) \Bigg\},  \mbox{ where, }    \\
u(\tau, x(\tau)) &=& 
\sup_{ a_k^i \in L^\infty [ \tau, T] }  J   (\tau, x(\tau), \  a_k^i ).
\end{eqnarray*}}
Observe that the function $J   ( \tau, x(\tau),   a_k^i)$ (see equation \eqref{eqn_opt_ctrl}) has   same structure as function  $J_h   (t, x,   a)$ in Lemma \ref{lemma 3}
with  $h = h^i_k$, $t = \tau$,  $x = x(\tau)$ and $a= a_k^i $ and continuity follows by Theorem \ref{thm_existance}. Hence we have,


\vspace{-2mm}
{\small \begin{eqnarray*}
u(\tau, x(\tau) ) & =& e^{-x(\tau)}
\left [ u(\tau, 0) - \nu x(\tau) \right  ], \mbox{ and so }
\\
%
u(t,0 ) &=&\sup_{ a_k^i \in L^\infty [t, \tau]}\Bigg\{\int_{t}^{\tau} \big (h^i_k(s)    - \nu   x (s)  \big )  a_k^i e^{-x(s)}ds  \\
&& \hspace{19mm}  + e^{-x(\tau)}[ u(\tau, 0 )-\nu x(\tau)]   \Bigg\}.
\end{eqnarray*}}
In the above equation  if we  consider zero policy, i.e., if we  consider
policy $a^i_k([t,\tau]) \equiv 0$ (basically  $a^i_k(s)= 0$ for all $s \in [t, \tau]$), then   clearly   $x(\tau)=0$, the first term (integral) in the above supremum is zero and  hence:
\TR{
$$
u(t,0)\ge u(\tau,0 ).  
$$}{$
u(t,0)\ge u(\tau,0 ).  
$ }
This 
  implies,
$
v^i_k (z^i_k;\pi^j)  \ge v^i_k (\bar{z}^i_k;\pi^j)  ,  
$ as  $ u(\tau, 0)$ is the value function of optimal control problem when the control starts in the state $\bar{z}^i_k$.  \eop

 \Cmnt{
\newpage
{\bf Proof of Lemma \ref{lemma_bang_bang_policy}:} For $M$  lock case, the best response  of agent $i$ against any given strategy of  player $j$, can be solved using (M stage) dynamic programming (DP) equations as below

 \vspace{-4mm}
{\footnotesize \begin{eqnarray}
v^i_{k} (\ z_k; \pi_{j} )\hspace{-1mm}  &\hspace{-1mm}  = \hspace{-1mm}  & \hspace{-1mm} 0 \mbox{ if }  k > M \mbox{ or if } \tau_{k-1}^i >  T \mbox{ or if }  z_k=(0, \tau_{k-1}),\label{eqn_dp_1} \\ \mbox{ and else,}\nonumber   \\
v^i_k (\ z_k; \pi_{j}) \hspace{-1mm}  &\hspace{-1mm}  = \hspace{-1mm}  & \hspace{-1mm} 
 \sup_{\ a_k}  \left \{ r_k^i (\ z_k, \ a_k; \pi_{j}) + E[  v^i_{k+1} (\ z_{k+1}; \pi_{j}) ) | \ z_k, \ a_k] \right \} \label{eqn_dp_2} . 
\end{eqnarray}}
We will show each stage of DP equation is optimized by a bang bang policy.

For any lock $k$, given any state $z_k$, the $k$-th stage DP equation can be re-written as below:

\begin{equation}\label{eqn_value_func_2}
    v_k^i(\ z_k; \pi_{j})  \ =  \sup_{\ a_k^i \in L^\infty }  J_{z_k} (\ a_k^i; \ \pi^j)   
\end{equation}
where the cost $J_{z_k}$ is defined as below:

\textcolor{red}{\footnotesize{
\begin{eqnarray}\label{eqn_cost}
    J_{z_k}(a_k^i; \pi^j) \hspace{-14mm}\\
    &=&
     \int_{\tau_{k-1}}^{T}     a_k^i ( t) e^{-\bar{a}_k^{i}(t)} \Bigg(  \sum_{l_{k+1}\in\{0,1\}}
    P^{i}(l_{k+1}/t,\pi_j) \nonumber \\ &&
   \hspace{36mm} \bigg [ c^i_k l_{k+1}  + 
    \ v_{k+1}^i( \underbrace{(l_{k+1}, t)}_{z_{k+1}},\pi^j) \bigg ]  \Bigg) dt  \nonumber \\
    &&-\nu \int_{\tau_{k-1}}^{T}\bar{a}_k^{i}(t)a_k^{i}(t) e^{-\bar{a}_k^{i}(t)}dt -\nu \bar{a}_k^{i}(T)e^{-\bar{a}_k^{i}(T)}, \nonumber 
\end{eqnarray}}
where
$\bar{a}^i_k  (s)   :=  \int_0^s a^i_k ( t )  dt =\bar{a}+ \int_{\tau_{k-1}}^s a^i_k ( t )  dt$ and $\bar{a}=0$
}
where $ P^{i}(l_{k+1}/t,\pi_j)$ is the probability that agent $i$ will be in state $z_{k+1}=(l_{k+1},t)$ at time $t$(the time of contact of $k$th lock) when other player uses $\pi^j$.
 In our case, $l_{k+1}=0$ means the other agent has contacted  before $t$, i.e., $z_{k+1}=(0,t)$. From equation (\ref{eqn_dp_1})  $\ v_{k+1}^i(z_{k+1},\pi^j)=0$, for such $z_{k+1}$, for any value of $\tau_{k-1}$ and hence, $$\ v_{k+1}^i((0,t),\pi^j)=0, \mbox{for any  }t$$ Basically, after the failure, the agent will stop trying, and value function will be zero there after. And $l_{k+1}=1$ means, the other agent has definitely not contacted the first lock before time $t$ (actually, agent $i$ found the lock first), and further, because we have $\tau_{k+1}$ approx $t$ with $t\leq T$. Thus for our case, $ P^{i}(l_{k+1}/t,\pi_j)$ is the probability that other agent has not contacted the lock before time $t$ for $k=1$,  and for all $k$ also. But the fact that for $k > 1$ and $z_k = (1, \tau_k)$ implies the other agent is not successful with first lock and hence, for any $k>1$, $$ P^{i}(l_{k+1}/t,\pi_j) = 1_{\{l_{k}= 1\}}$$ 
In view  of the above, equation (\ref{eqn_cost}) modifies to

\textcolor{red}{\begin{eqnarray}
       J_{z_k}(a^i; \pi^j)=\nonumber \\
       && \hspace{-20mm}\int_{\tau_{k-1}}^{T}  a_k^i ( t)e^{-\bar{a}_k^{i}(t)}  \bigg(\Big( c^i_k+ v_{k+1}^i\big((1,t),\pi^j\big)\Big)
       P_{f}^{j}(t/\pi_j)\nonumber\\
       &&\hspace{10mm}
       -\nu \bar{a}_k^i(t)\bigg) dt -\nu \bar{a}_k^{i}(T) e^{-\bar{a}_k^{i}(T)},
\end{eqnarray}}
\Cmnt{
The above problem resembles to a finite horizon optimal control problem, where $J$ can be rewritten as
\begin{eqnarray}
\label{Eqn_opt_cntrl_cost}
J_h(t,x,a) =\int_{t}^{T} \left ( h^j (s)    - \nu   x (s)  \right ) a_k^i (s)\exp(-x(s))ds   +g(x(T)),
\end{eqnarray}
with state process
$$
\stackrel{\bm\cdot}{x}  (s) =    a_k^i ( s)  \mbox{ and thus }  x(s) = \int_0^s a_k^i (l) dl=\Tilde{x}+\int_t^s a_k^i (l) dl = \bar{a}_k^i (s)\ ;\ \mbox{\textcolor{red}{$\Tilde{x}=
0, x(t)=x$}}
$$a given function 
$$
h^j(s) := \ v_{k+1}^i(1,s,\pi^j)P_{f}^{j}(s/\pi_j),
$$and
 with terminal cost 
\begin{equation*}
    g(x) = - \nu x\exp (-x ). \nonumber
\end{equation*}}

\textcolor{red}{  One can write the  optimization problem of the above $k$-th stage DP equation as the value function 
of 
an appropriate optimal control problem  with details  as follows,
\begin{eqnarray*}
v_k^i(\ z_k; \pi_{j})  &=&  u(t, 0)
 \\
u(l, x) &:= & \sup_{a_k^i \in L^\infty [l, T] } J (l, x, \ a_k^i )  \mbox{ where, }  \\
 J   (l, x, \ a_k^i) & :=&\int_{l}^{T} \left ( h_k^j (s)    - \nu   x (s)  \right ) a_k^i ( t) e^{-x(s)}ds   +g(x(T )),  \mbox{ with, }\\ 
h_k^j(s) &:=&  (c^i_k + v_{k+1}^i( (1,s); \pi^j) )P_{f}^{j}(s/\pi_j),
\end{eqnarray*}
with state process (for any $s \in [l, T]$)
$$
\stackrel{\bm\cdot}{x}  (s) =    a_k^i ( s), \mbox{ with } x(l)= x  \mbox{ and thus }  x(s) = x+ \int_l^s a_k^i (l) dl 
$$and
 with terminal cost 
\begin{equation*}
    g(x) = - \nu x e^{-x}. \nonumber
\end{equation*}
 }
The equation (\ref{eqn_value_func_2}) can be rewritten as

$$u(t,0)= \sup_{ \ a (.) }   \left \{  \int_t^T  L(s, \ x (s) ,  \ a(s) )  ds   + g(x (T)  ) \right \},$$ { with }  \\
$$L(s, \ x, \ a) = 
  \left ( h^j (s)    - \nu   x  \right ) a e^{-x},$$   and   $$g(\ x) = - \nu  x e^{- x}.$$
The optimal control problems can be solved by solving respective HJB PDEs \cite{feller}.  
Thus, we need the solution of the following (Hamiltonian Jacobi) HJB PDE,
$$
u_s  (t, 0) + \max_{a \in [0, \beta]}   \ a\{  ( h^j (s)    - \nu   x   ) e^{-x}    +  u_x (t,0) \}  = 0,
$$
with boundary condition 
$$
u(T, x) =  g( x) =  -\nu \ x e^{- x}.
$$
where, $u(t,0)$ represents the value function, and  $u_s$, $u_x$ are its partial derivatives.
The PDE\footnote{There are more technicalities, such that the solution of PDE should exist.} shows that the maximizer will be at either of the extreme points. 
%
%
So, the best response against any policy of opponent, i.e, the acceleration process $a_{k^*}^i$ to acquire lock $k$ against any policy $\pi^j$ of opponent will be a bang bang policy for every possible $z_k$, and for every $k$.

\Cmnt{
\textbf{Proof of lemma \ref{Lemma_value_function_mono}}

 For any lock $k$, given any state $z_k$, the $k$-th stage DP equation can be re-written as below:
\begin{equation*}
    v_k^i(\ z_k; \pi_{j})  \ =  \sup_{\ a_k^i \in L^\infty } 
      J_{z_k} (  \ a_k^i; \ \pi^j)  , 
\end{equation*}
where $J_{z_k}$ is defined as in equation (\ref{cost}).   One can write the  optimization problem of the above $k$-th stage DP equation as the value function 
of 
an appropriate optimal control problem (see equation \ref{eqn_opt_ctrl}) with details  as follows,
\begin{eqnarray*}
v_k^i(\ z_k; \pi_{j})  &=& u(t, 0; z_k) \mbox{ where  the optimal control problem is defined by } \\
u(l, x;z_k) &:= & \sup_{a_k^i \in L^\infty [l, T] } J (l, x, \ a_k^i; \z_k )  \mbox{ where, }  \\
 J   (l, x, \ a_k^i; \ z_k) & :=&\int_{l}^{T} \left ( h_k^j (s)    - \nu   x (s)  \right ) a_k^i ( t)e^{-x(s)}ds   +g(x(T )),  \mbox{ with, }\\ 
h_k^j(s) &:=& \hspace{-3mm} (c^i_k + v_{k+1}^i( (1,s); \pi^j) ) \left  (1_{ k > 1} + 1_{k=1}P_{f}^{j}(s/\pi_j) \right ),  \ \ z_{k+1,s} := \mbox{   $(1,s)$  }
\end{eqnarray*}
with state process (for any $s \in [l, T]$)
$$
\stackrel{\bm\cdot}{x}  (s) =    a_k^i ( s), \mbox{ with } x(l)= x  \mbox{ and thus }  x(s) = x+ \int_l^s a_k^i (l) dl 
$$and
 with terminal cost 
\begin{equation*}
    g(x) = - \nu x e^{-x}. \nonumber
\end{equation*}
  
  From the Dynamic programming principle \cite[Theorem 5.1?]{feller}, we can write it as follows;
 \begin{eqnarray}
u(t,0; R_{z_k} )  &=&\sup_{ a_k^i \in L^\infty [t, \tau] }\Bigg\{ \int_{t}^{\tau} \left ( h_k^j (s)    - \nu   x (s)  \right )  a_k^i ( s)e^{-x(s)}ds   + u(\tau,x(\tau) ;R_{z_k}) \Bigg\},  \mbox{ where, }  \nonumber  \\
u(\tau, x(\tau); R_{z_k} ) &=& 
\sup_{ a_k^i \in L^\infty [ \tau, T] }  J   (\tau, x(\tau), \  a_k^i; \ R_{z_k} )
\end{eqnarray}

Observe that the function $J   (\tau, x(\tau), \  a_k^i; \ R_{z_k} )$ has the same structure as in Lemma \ref{lemma 3}
with  $h = h^j_k$, $t = \tau$,  $x = x(\tau)$ and $a= a_k^i $.

In a similar way, when we start with $z'_k$ the resulting optimal control problem again has the almost same structure as the optimal control problem of Lemma \ref{lemma 3}, i.e., we again have $h = h^j_k$, $t = \tau$ and $a= a_k^i$, but now $x = 0$.    By Lemma \ref{lemma 3}, we have that 
$$
u(\tau, x(\tau); R_{z_k} )  = e^{-x(\tau)}
\left [ u(\tau, 0; R_{z_k} ) - \nu x(\tau) \right  ]
$$
Note that $u(\tau, 0; R_{z_k} ) = u(\tau, 0; R_{z'_k} )$ by hypothesis.
 Thus, 
 \begin{eqnarray}
u(t,0; R_{z_k} ) =\sup_{ a_k^i \in L^\infty [t, \tau]}\Bigg\{\int_{t}^{\tau} \left (h^j_k(s)    - \nu   x (s)  \right )  a_k^i e^{-x(s)}ds   + e^{-x(\tau)} u(\tau, 0; R_{z'_k} )-\nu x(\tau)] ) \Bigg\}, 
\end{eqnarray}
In particular consider zero policy, i.e., consider
policy $a^0([t,\tau]) \equiv 0$, where, $a^i_k(s)= 0$ for all $s \in [t, \tau]$. Under this policy clearly,  we have $x(\tau)=0$, the first integral in the above supremum is zero and  hence we have
$$
u(t,0; R_{z_k} ) \ge J(t, 0, a^0) = u(\tau,0; R_{z'_k} ).  
$$
which implies,
$$
v_k (z_k;\pi^j)  \ge v_k (z'_k;\pi^j)  
$$
}}
\Cmnt{
@@@@@@@@@\\
We have a unique optimizer for each $k$ as a function $t_{k-1}$ as given below, which are defined backward recursively:
\begin{eqnarray}
\theta_k^* (t_{k-1} )  &=& 
\left \{ 
\begin{array}{cccc}
   t_{k-1}      &  \mbox{ if }       t_{k-1}    >   \theta_k^*     \\
\theta_k^* &  \mbox{ else,   where  } 
   \end{array} \right . \nonumber  \\
 \theta_k^*  &:=&     \inf_{t \ge 0 }   \{   c^i_k+\Upsilon_{k+1}^* (t)        \le   \nu  \}   \wedge T.
\end{eqnarray}
The cost to go starting from $t_{k-1}$ is strictly decreasing for all $t_{k-1}  < \theta^*_k$, after which it remains at 0: 

\vspace{-4mm}
{\small\begin{eqnarray*}
\Upsilon_M^{*}(t_{M-1}) &=& 1_{c_M \ge \nu} e^{ \beta t_{M-1} } \ \int_{ t_{M-1}}^{ T} (c_M - \nu )\beta e^{-\beta  s  }d s \\
&&
\mbox { and for }  k < M \\
\Upsilon_{k}^* (t_{k-1})   &=& \int_{t_{k-1}}^{  \theta^*_k} \left [  c^i_k+\Upsilon_{k+1}^* (t_k )  
   - \nu      \right ]
\beta e^{-\beta t_k  }d t_k
\end{eqnarray*}}
@@@@@@@@@
\begin{thm}
Against any fixed strategy of the opponents, there exists a unique threshold strategy in the best response, i.e., for every lock $k$, there exists a unique threshold policy which is in the best response against any fixed strategy of the opponent. In other words, for every $k$, there exists a unique time threshold as an optimizer.
\end{thm}
{\bf Proof of the theorem:} 
To prove the uniqueness of the threshold policy, we will prove that there is a unique time threshold which is optimal, and the uniqueness of time threshold will imply the uniqueness of threshold policy. Note that the optimal time threshold for any lock $k$ is independent of the contact epoch $\tau_{k-1}$ of the previous lock,  and is given by:
\begin{eqnarray*}
\theta^{i*}_k&\in&  \arg \max_{\theta_k} \int_{0}^{\theta_k} (c^i_k+\Upsilon_{k+1}^* (t)-\nu)\beta e^{-\beta t}dt.\\
\end{eqnarray*}
To prove the uniqueness of  optimizer $\theta^{i*}_k$, we will first prove the function $\Upsilon_{k+1}^* (t) $    is a strict decreasing function of $t$.  
Infact, we will prove the function  $\Upsilon_{k}^* (t) $   is strictly decreasing for every $k$. We will prove this using induction, for $k=M$ we have, 	
\begin{eqnarray*}
\Upsilon_M^{*}(t) &=& \int_{ 0}^{\theta^{i*}_M-t} (c_M - \nu )\beta e^{-\beta s }d s \\
\end{eqnarray*}
where $\theta^{i*}_M=T$ if ${c_M >\nu}$ and zero otherwise. Observe that this is a strictly decreasing function of $t$ till $\Upsilon_M^{*}(t)$ becomes zero and it remains zero afterwards. Assume the function $\Upsilon_k^{*}(t)$ is strictly decreasing function till  is true for $i=M, M-1\dots k+1$ and now, if we prove this holds true for $i=k$ we are done, we will have the desired result. For $i=k$, we have
\begin{eqnarray*}
\Upsilon_{k}^* (t) = \left \{ 
\int_{0}^{ \theta^*_k- t } \left [  c^i_k+\Upsilon_{k+1}^* (t+s)    
   - \nu      \right ]
\beta e^{-\beta s}d s \right    \},
\end{eqnarray*}
By the induction hypothesis, $\Upsilon_{k+1}^* (t) $ is a strict decreasing function of $t$ which implies $\Upsilon_{k+1}^* (t+s) $ is also a strict decreasing function of $t$. We have the integrand is a decreasing  function of $t$, and so is the integral, it implies $\Upsilon_{k}^* (t)$ as a decreasing function of $t$.
And hence, we  have  $\Upsilon_{k}^* (t) $ for every k is a strict decreasing function of $t$ till it becomes zero.\par 

Now we will prove the uniqueness of the optimizers (optimal time thresholds) as follows: we have for any k the optimal threshold is given by  $\theta^{i*}_k$
\begin{eqnarray*}
\theta^{i*}_k&\in&  \arg \max_{\theta_k} \int_{0}^{\theta_k} (c^i_k+\Upsilon_{k+1}^* (t)-\nu)\beta e^{-\beta t}dt.\\
\end{eqnarray*}
On maximizing the above equation with respect to $\theta_k$ we have;
$$\theta_k^*  =    \inf_{t \ge 0 }   \{   c^i_k+\Upsilon_{k+1}^* (t)        \le   \nu  \}   \wedge T.$$
where $\Upsilon_{k+1}^* (t)$ is a strict decreasing fuction of $t$ till it becomes zero, which implies $\theta_k^*$ is unique. So, for every $k\in \{2,3,\dots,M\}$ we have the unique optimizer. For k=1, we have 
\begin{eqnarray*}
\theta^{i*}_1
&\in &  \arg \max_ {\theta_k} \int_{0}^{\theta_1} \left [   ( c_1 + \Upsilon_2^* (t)  ) P(\Psi > t)    
   - \nu      \right ]
\beta e^{-\beta t_1}dt  
\end{eqnarray*}
in the above equation the integrand is a strictly  decreasing function of $t$, which implies it is a strict concave function which implies the uniqueness of $\theta^{i*}_1$.
Note that the control for the $k-th$ lock starts at $\tau_{k-1}$ and if $\tau_{k-1}$ is greater than $\theta^{i*}_k$, it means the agent will not try for the $k-th$ lock and will stop at $k-1$th lock. Hence, the optimal threshold can be written as a function of $\tau_{k-1}$ as follows:
\begin{eqnarray*}
\theta_k^* (t_{k-1} )  \in \arg\max_{\theta_k} \left \{ 
\int_{0}^{\theta_k - t_{k-1} } \left [  c^i_k+\Upsilon_{k+1}^* (t_k+t_{k-1})    
   - \nu      \right ]
\beta e^{-\beta t_k}dt_k \right    \},
	\end{eqnarray*}}

\begin{lemma}\label{lemma 3}
\label{Lemma_ic_change}
Let  $J_h(t,x,a )$   be a function of the form

{\small
 \begin{eqnarray}\nonumber
\label{Eqn_opt_cntrl_one_lock}
J_h(t,x,a  ) =\int_{t}^{T} \left ( h (s)    - \nu   x (s)  \right ) a ( s) e^{-x(s)}ds   
- \nu x(T) e^{-x(T) }, \nonumber
\end{eqnarray}}defined using  continuous function $h(\cdot)$ and  
 state process
$$
\stackrel{\bm\cdot}{x}  (s) =    a  ( s), \mbox{ with initial condition, }  x(t)=x
.
$$
Define $u(t,x) := \sup_{\ a \in L^\infty }  J_h (t,x,a)$, then  we have:
\begin{enumerate}[(i)]
    \item \vcmnt{ $J_h(t,x,a )=e^{-x}[J_h(t,0,a )-\nu x]$}
    \item $u(t,x)=e^{-x}[u(t,0)-\nu x]$
    \item The optimal policy  $a^*(\cdot) $ is independent of $x$.
\end{enumerate}                                                                                                                                                                                                                 
    \end{lemma}{}
\noindent{\bf Proof:} 
By  change of variables
$x(s)=x+\Tilde{x}(s),$  
$$
\stackrel{\bm\cdot}{\Tilde{x}}  (s) =    a ( s)  \mbox{, i.e.,  }  \Tilde{x}(s) = \int_t^s a ({\tilde s}) d{\tilde s}  , \mbox{ and }  \Tilde{x}(t)=0,
$$we get,

\vspace{-6mm}
{\small
  \begin{multline*} 
    J_h(t,x,a )  =\int_{t}^{T} \Big( h (s)    - \nu  \big ( x+\Tilde{x} (s) \big)  \Big ) a ( s) e^{ -(x+\Tilde{x}(s))} ds\\
     -\nu  \big(x+\Tilde{x}(T)\big) e^{-x-\Tilde{x}(T)} \\ \TR{
     =e^{-x} \Bigg ( \int_{t}^{T} \Big( h (s)    - \nu   \big( x+\Tilde{x} (s)\big )  \Big ) a ( s) e^{ - \Tilde{x}(s)} ds\\
     -\nu  \big(x+\Tilde{x}(T)\big) e^{-\Tilde{x}(T)} \Bigg ) \\  }{}
     = e^{-x} \Bigg ( \int_{t}^{T} \big( h (s)    - \nu  \Tilde{x} (s)   \big ) a ( s)e^{ - \Tilde{x}(s) }ds - \nu \Tilde{x}(T)  e^{-\Tilde{x}(T)} \\
    \hspace{10mm} - \nu x \int_{t}^{T}    a( s)e^{ - \Tilde{x}(s)} ds
     -\nu  x e^{-\Tilde{x}(T)} \Bigg ) \\
     = e^{-x} \left ( J_h(t, 0,a) - \nu x \right ). \hspace{45mm} 
 \end{multline*}}
 The last equality  follows because the sum of  two terms (probabilities) is 1.   This completes part~(i).
 For  part (ii), from the above  equation we have,
   $$
    u(t,x)=\sup_{a\in L^{\infty} } J_h(t, x,a) = \sup_{a\in L^{\infty}}   e^{-x} \left [  J_h (t, 0,a) - \nu x \right  ]
   $$
 and hence   we have,
$
u(t,x)=e^{-x}[u(t,0)-\nu x].
$
This proves part (ii). 
Further it is clear from above that the optimal policy $a^*(\cdot)$ remains the same for all initial conditions $x$ and this proves part (iii).  \eop\par

\TR{
\section*{Appendix R: Reduced Game}
{\bf Proof of Theorem \ref{Thm_Upsilon_two}:} 
Define the following,   for any $2 \le l \le M$ and any MT-strategy $\pi =(\theta_l \dots\theta_M)$ (see   \eqref{Eqn_Strategy}):
$$
\gamma_l (\theta_l \cdots, \theta_M;  t ) :=   \sum_{k = l}^M E[ r_k^i(\ z_k^i, \theta_k )   |  z^i_l  = (1, t) ] .
$$Note these objective functions (with $l \ge 2$) do not depend upon the strategies of opponents. Define
$$
\Upsilon_l^* (t) := \sup_{\theta_l, \cdots, \theta_M} \gamma_l (\theta_l \cdots, \theta_M;  t ) ,
$$and observe that
$$
\Upsilon_2^{i*} (t)= \Upsilon_2^* (t)= \sup_{\theta_2, \cdots, \theta_M} \gamma_2 (\theta_2 \cdots, \theta_M;  t ) .
$$
Further, by definition,   for any given $(\theta_{k-1}, \cdots, \theta_M)$ we have:

\vspace{-4mm}
{\small \begin{eqnarray}
 \gamma_{k-1}  (\theta_{k-1} \cdots, \theta_M;  t ) \hspace{-25mm} && \label{Eqn_gamma_k} \\
 &=&  \int_t^{\theta_{k} }  \left (    c^i_{k-1} - \nu  +    \gamma_{k} (\theta_{k} \cdots, \theta_M;  s  )   \right )      \beta^i  e^{-\beta^i (s-t) }  ds \nonumber \\
 &\le & \hspace{-1mm}   \int_t^{\theta_k } \hspace{-1.mm}  \left (    c^i_{k-1} - \nu  +   \sup_{\theta'_{k} \cdots, \theta'_M} \gamma_{k}  (\theta'_{k} \cdots, \theta'_M;  s  )   \right )      \beta^i  e^{-\beta^i (s-t) }  ds, \nonumber 
 \end{eqnarray}}with $\gamma_{M+1}\equiv 0$. Observe that $ \sup_{\theta'_k \cdots, \theta'_M} \gamma_k  (\theta'_k \cdots, \theta'_M;  s  )  $ is the problem of finding the best response against silent opponent   (i.e., when none of the  opponents are attempting) with $M-k$ locks.   By applying  Theorem \ref{thm_monotonocity_of_thrshold}  one can choose the MT-strategy as a best response, i.e.,
  $$ \sup_ {\pi= \{ a_k (z_k) \} \cdots,  \{ \a_M (z_M) \} } \gamma_k (a_k \cdots, a_k;  s  ) =  \max_{\theta'_k \cdots, \theta'_k} \gamma_k (\theta'_k \cdots, \theta'_M;  s  ),$$
 also optimal  $(\theta^*_k \cdots, \theta^*_M)$   do not depend  upon $s$, i.e., 
%

\vspace{-4mm}
{\small
\begin{eqnarray}
\Upsilon_k^* (s) : =  \max_{\theta'_k \cdots, \theta'_M} \gamma_k (\theta'_k \cdots, \theta'_M;  s  ) =    \gamma_k (\theta^*_k \cdots, \theta^*_M;  s  ) \mbox{ for all }  s. 
\label{Eqn_Upsilon_star_k}
\end{eqnarray}}
And hence we have,
\begin{eqnarray}
\gamma_{k-1}  (\theta_{k-1} \cdots, \theta_M;  t ) \hspace{-27mm}  \nonumber\\
 &\le &    \int_t^{\theta_{k-1} }  \left (    c^i_{k-1} - \nu  + \Upsilon_k^* (s)\right )      \beta^i  e^{-\beta^i (s-t) }  ds. \nonumber\\ 
 &\le &  \sup_ {\theta'_{k-1}}   \int_t^{\theta_{k-1} }  \left (    c^i_{k-1} - \nu  + \Upsilon_k^* (s)\right )      \beta^i  e^{-\beta^i (s-t) }  ds.  \hspace{4mm}
   \label{Eqn_gamma_k}
 \end{eqnarray}
Consider the case with $k = M$. In this case  clearly 
\begin{eqnarray*}
\Upsilon_M^* (s)  &: =&  \max_{ \theta'_M} \gamma_M (  \theta'_M;  s  ) \\
&=& \max_{ \theta'_M}    \int_t^{\theta'_M }  \left (    c^i_M - \nu  \right )      \beta^i  e^{-\beta^i (s-t) }  ds.
\end{eqnarray*}
 The integrand is a strictly decreasing function, which implies the integral is strictly concave and hence has a unique maxima $\theta^*_M$, as given below along with optimal $\Upsilon_M^*$:
\begin{eqnarray*}
\Upsilon_M^* (s)&=& \gamma_M   (\theta^*_M)  \mbox{ with }\theta^*_M  =  T  \indc{ c_M > \nu} \\
&=&      (c^i_M-\nu) \left ( 1 -  e^{- \beta^i (T-s) }  \right ) \indc{ c_M > \nu}  \mbox{ for all }  s. 
\end{eqnarray*}

 Observe that the function  $\Upsilon_M^*  (\cdot)$ is a strictly decreasing for $s<\theta^*_M$ or remains at zero for all $t$ if  $\theta^*_M = 0$ and further the coefficient $\theta^*_M$ 
 is unique.
In other words,  the function is strict decreasing for all $s \le \theta^*_M$ and remains at 0 after  (unique) $\theta^*_M$.
Assume the same  holds true for all  $k=M, M-1 \dots p+1$ (for any $p$) and then consider  $k=p$.  From equation \eqref{Eqn_gamma_k} 
\begin{eqnarray*}
 \gamma_{p} (  \theta'_{p},\dots  \theta'_M; t  )\hspace{-3mm}  & \le & \hspace{-2mm}  \sup_ {\theta'_{p}}    \int_t^{\theta'_{p} }  \left (    c^i_{p} - \nu  +   \Upsilon^*_{p+1} (   s  )   \right ) \beta^i  e^{-\beta^i (s-t) }  ds.
\end{eqnarray*}
Since $\Upsilon^*_{p+1}(\cdot)$ is non-increasing function, the  integrand in the above inequality  is strictly decreasing with $s$. This implies the upper bounding integral is strictly  concave and hence has an unique maximizer, in fact the maximizer equals:
 $$  \theta_{p}^{*} =    \inf   \{  t \ge 0 :  c^i_{p}+\Upsilon_{p+1}^{*} (t)        \le   \nu  \}   \wedge T, \mbox{ with } \inf   \emptyset  := 0.
$$
Thus we have for any  $ \theta'_{p}, \dots \theta'_M$:
\begin{eqnarray*}
\gamma_{p} (  \theta'_{p},\dots  \theta'_M;  t  ) &  \le & \gamma_{p} (  \theta^*_{p}, \dots \theta^*_M;  t  )  \mbox{ and hence }
\end{eqnarray*}
 \begin{eqnarray*}
\Upsilon_{p}^* (t)  & = & \gamma_{p} (  \theta^*_{p},\dots  \theta^*_M;  t )  \\
 &=&     \int_t^{\theta^*_{p} }  \left (    c^i_{p} - \nu  +   \Upsilon^*_{p+1} (   s  )   \right ) \beta^i  e^{-\beta^i (s-t) }  ds
\end{eqnarray*}
From the above it is clear that the function $\Upsilon_{p}^* (t) $   is strictly decreasing with $t$ for all $t \le \theta^*_{p} $  and 
$\Upsilon_{p}^* (t)  = 0$ for all $t > \theta^*_{p} .$  The proof is complete by backward induction, with ${\bar \Upsilon}_k^{i*} = \Upsilon_k^* $ 
and $\theta^{i*}_k = \theta^*_k$ for all $k$.
 \eop
}{}

\Cmnt{
\section{Numerical Simulations}


%
%
\Cmnt{
Consider the symmetric case, i.e.,  $\beta^i=\beta$  and $c^i_k=c_k\mbox{ for all } i\in N$. The uniqueness of NE given by Theorem \ref{thm_unique_NE} implies for any $k$, $\theta_k^i=\theta_k^j=\theta_k^*\mbox{ for all } i,j\in N$.\par
 The NE for $1$ lock and $n$ player symmetric game using Theorem  \ref{thm_unique_NE} is given by time thresholds $(\theta^*\dots\theta^*)$  with $\theta^*$ defined as:
\begin{equation*}
\theta^* := \left \{ \begin{array}{lllll}
   -\frac{  \ln ( \frac{\nu}{c_1} ) } { (n-1) \beta }   & \mbox{ if }   e^{-\beta (n-1) T}  \le  \nu\\
   T   &\mbox{ else.}
\end{array}      \right .
\end{equation*}
which also proves the conjecture given in \cite{Mayank} with $c_1=1$.
}

%
%
%
%


}

\Cmnt{
\subsection{
When others are silent}
In this case, $ P(\Psi > t) = 1$ for all $t$ and hence one needs to maximize 

{\small \begin{eqnarray*}
\Upsilon_{k}^{i*} (t_{k-1}) = \max_{\theta_k} \left \{ 
\int_{t_{k-1}}^{\theta_k} \left [  \Upsilon_{k+1}^{i*} (t_k)    
   - \nu      \right ]
\beta e^{-\beta t_k}dt_k \right    \}
\end{eqnarray*}}
the maximizer would be:
$$
\theta_k^{i*} (t_{k-1} ) = t_{k-1} \wedge \theta^{i*}_k
$$

 Since 
$\Upsilon^{i*}  $ is decreasing function, it is clear that
$$
\theta_k^{i*} = \sup_t \left   \{ \Upsilon_{k+1}^{i*} (t)    
   >  \nu \right  \}
$$
When $k = M$ (trying the last lock) then clear that $ \Upsilon_M^{i*} (t)   = c^i_M$ for all $t$ and hence 
$$
\theta_M^{i*} =  T 1_{ c^i_M > \nu }.
$$
And then 
$$
\phi^{i*}_M = (c^i_M - \nu )  (1 - e^{ - \beta^i T  } )   1_{ c^i_M > \nu }.
$$
When $k = (M-1)$, then (from second silence Theorem  or from above statement) we have
$$ \Upsilon_{M-1}^{i*} (t)   = c_{M-1} +  (c^i_M - \nu )  (1 - e^{ - \beta^i (T-t)  } )   1_{ c^i_M > \nu } \mbox{ for all } t$$ 
which is clearly decreasing with $t$.
Thus 
$$
(1 - e^{ - \beta^i (T-\theta_{M-1}^{i*} )  } )   =  \frac{  \nu -  c_{M-1}  } {c^i_M - \nu }
\mbox{ or }  $$
$$
\theta_{M-1}^{i*}  =  \frac{1}{\beta^i} \left ( T -  \log \left (  \frac{c^i_M - \nu } {  c^i_M - 2\nu +  c_{M-1}  } \right ) \right )
$$
We now have
\begin{eqnarray*}
 \Upsilon^{i*} (t)  & = & c_{M-2} + \sum_{l=0}^1  (c_{M -l}- \nu )  (1 - e^{ - \beta^i (\theta_{M-1}^{i*}  -t)  } )  \\ 
 &&   
-  (c^i_M - \nu )  e^{-\beta^i  T } ( \theta_{M-1}^{i*}  - t)
  \mbox{ for all } t \le \theta_{M-1}^{i*}  
\end{eqnarray*}
(Assuming some kind of monotonicity again)
Now, $\theta_{M-2}^{i*}$ satisfies

 {\small \begin{eqnarray*}
 \sum_{l=0}^1  (c_{M -l}- \nu )   e^{ - \beta^i (\theta_{M-1}^{i*}  - \theta_{M-2}^{i*} )  }   
 +  (c^i_M - \nu )  e^{-\beta^i  T } ( \theta_{M-1}^{i*}  - \theta_{M-2}^{i*} )  \\  =   \sum_{l=0}^2  (c_{M -l}-\nu)
\end{eqnarray*}}
And 

{\small
\begin{eqnarray*}
 \Upsilon^{i*} (t)  & = & c_{M-3} + \left (  \sum_{l=0}^2  (c_{M -l}- \nu )   -  (c^i_M - \nu )  e^{-\beta^i  T }  \theta_{M-1}^{i*} \right )  (1 - e^{ - \beta^i (\theta_{M-2}^{i*}  -t)  } )  \\ 
 &&   
-  \sum_{l=0}^1  (c_{M -l} -\nu )  e^{-\beta^i  \theta_{M-1}^{i*} } ( \theta_{M-2}^{i*}  - t)   \\
&& +   (c^i_M - \nu )  e^{-\beta^i  T }  \frac{ ( \theta_{M-2}^{i*}  - t)^2  }{2 }
  \mbox{ for all } t \le \theta_{M-2}^{i*}   
\end{eqnarray*}}

\newpage

\vspace{-4mm}
{\footnotesize{\begin{eqnarray}
  P_f &=& \int_{0}^{\theta_1}\Bigg(\int_{t_1}^{\theta_2}\Big(\int_{t_2}^{\theta_3}\dots\big(\int_{t_{M-1}}^{\theta_M}\beta e^{-\beta (t_M-t_{M-1})}dt_M\big)\nonumber\\
 && \hspace{10mm}\dots\beta e^{-\beta (t_3-t_2)}dt_3\Big)\beta e^{-\beta (t_2-t_1)}dt_2\Bigg)\nonumber\\
 &&\hspace{22mm} e^{-\beta_{j}(t_1\wedge\psi_1)}\beta e^{-\beta t_1}dt_1\nonumber
\end{eqnarray}}

 {\begin{eqnarray}
  P_f &=& \beta^M \int_{0}^{\theta_1} \int_{t_1}^{\theta_2} \int_{t_2}^{\theta_3}\dots \int_{t_{M-1}}^{\theta_M}  e^{-\beta t_M}dt_M  dt_{M-1}  \cdots dt_2 \nonumber\\ 
 &&\hspace{32mm} e^{-\beta_{j}(t_1\wedge\psi_1)}  dt_1\nonumber
\end{eqnarray}} }

{\color{blue}{\small
\begin{eqnarray*}
&=&\beta\int_0^{\theta_1} \Big(\beta^{M-2}
\int_{t_1}^{\theta_2} \int_{t_2}^{\theta_3}\dots \int_{t_{M-2}}^{\theta_{M-1}}\\
&& \hspace{20mm} \left ( e^{-\beta \theta_M }  - e^{-\beta t_{M-1}} \right )dt_{M-1}\\
&&\hspace{25mm}  \cdots dt_2\Big) e^{-\beta_{j}(t_1\wedge\psi_1)}  dt_1\\
&=&\beta\int_0^{\theta_1} \Big( \sum_{k=0}^{M-2} \mbox{\color{red}$($}(-1)^k \mbox{\color{red}$-\beta \nu \theta_{M-k-1})$} \beta^{M-k-2} e^{-\beta \theta_{M-k} }  \\
&& \hspace{10mm}
\int_{t_1}^{\theta_2}\dots \int_{t_{M-k-2} }^{\theta_{M-k-1}}   dt_{M-k-1} \cdots dt_2 \\
&&\hspace{15mm} + e^{-\beta t_1} (-1)^{M}\Big) e^{-\beta_{j}(t_1\wedge\psi_1)}  dt_1
\end{eqnarray*}
 } }
{\footnotesize
\begin{eqnarray*}
%
&=&\beta\int_0^{\theta_1} \Big( \sum_{k=0}^{M-2} (-1)^k \beta^{M-k-2} e^{-\beta \theta_{M-k} }  \\
&& \hspace{10mm}
\int_{t_1}^{\theta_2}\dots \int_{t_{M-k-2} }^{\theta_{M-k-1}}   dt_{M-k-1} \cdots dt_2 \\
&&\hspace{15mm} + e^{-\beta t_1} (-1)^{M}\Big) e^{-\beta_{j}(t_1\wedge\psi_1)}  dt_1
\end{eqnarray*}
  }
{\color{green}  
  All we need is the Lebegue measure/volume of the simplex 
  $$
 \left  \{  x_{i} \le \theta_i \forall 2 \le  i \le M-1, \ \ x_i \le x_{i+1} \forall 2\le  i \le M-2,   \ t_1 \le  x_2  \right \}  
  $$
  
  for M=4
\begin{eqnarray*}
  \{ t_1 \le x_2 \le \theta_2 , \  x_2 \le x_3 \le \theta_3  \} \hspace{-40mm} \\
 & =&  \{ t_1 \le x_2 \le \theta_2 , \  0  \le x_3 \le \theta_3  \}
  - \{ t_1 \le x_2 \le \theta_2 , \  0  \le x_3 \le x_2 \} \\
   & =& [t_1, \theta_2] \times [0, \theta_3] 
  -  1/2  [t_1,  \theta_2 ]^2   \\
\end{eqnarray*}

    $$
  \{ t_1 \le x_2 \le \theta_2 , \  x_2 \le x_3 \le \theta_3 \   \  x_3 \le x_4 \le \theta_4   \}
  $$
  }

and the expected cost will be


{\small{\begin{eqnarray}
 E[cost] &=& \beta^{M+1}\int_{0}^{\theta_1}\int_{t_1}^{\theta_2}\int_{t_2}^{\theta_3}\dots\int_{t_{M-1}}^{\theta_M} t_M e^{-\beta t_M} {dt_M}\nonumber\\
 &&\hspace{20mm}\dots dt_3{dt_2}e^{-\beta_{j}(t_1\wedge\psi_1)} dt_1\nonumber\\
 &&+ \beta^{M}\int_{0}^{\theta_1}\int_{t_1}^{\theta_2}\int_{t_2}^{\theta_3}\dots\int_{t_{M-2}}^{\theta_{M-1}} \theta_M e^{-\beta \theta_M} {dt_{M-1}}\nonumber\\
 &&\hspace{20mm}\dots dt_3{dt_2}e^{-\beta_{j}(t_1\wedge\psi_1)} dt_1\nonumber\\
&&+ \beta^{M-1}\int_{0}^{\theta_1}\int_{t_1}^{\theta_2}\int_{t_2}^{\theta_3}\dots\int_{t_{M-3}}^{\theta_{M-2}} \theta_{M-1} e^{-\beta \theta_{M-1}} \nonumber\\
 &&\hspace{20mm}{dt_{M-1}}\dots dt_3{dt_2}e^{-\beta_{j}(t_1\wedge\psi_1)} dt_1\nonumber\\
 &&\hspace{22mm}e^{-\beta_{j}(t_1\wedge\psi_1)}\beta dt_1\nonumber\\
 &&\hspace{15mm}\vdots\nonumber\\
  &&\hspace{15mm}\vdots\nonumber\\
  && + \beta^2 \int_{0}^{\theta_1}    \theta_2 e^{-\beta\theta_2}   
   e^{-\beta_{j}(t_1\wedge \psi_1)} dt_1\nonumber\\
  &&+\beta\theta_1 e^{-\beta \theta_1}\nonumber
\end{eqnarray}} }

{
\footnotesize
\begin{eqnarray*}
E[cost]&=& \sum_{k=1}^{M} \beta^{k} \theta_{k} e^{-\beta\theta_{k}}\int_{0}^{\theta_1}\int_{t_1}^{\theta_2}\\
&&\hspace{18mm}\dots\int_{t_{k-2}}^{\theta_{k-1}}dt_{k-1}\dots dt_2  e^{-\beta_{j}(t_1\wedge \psi_1)} dt_1\\
&&+ \beta^{M+1}\int_{0}^{\theta_1}\int_{t_1}^{\theta_2}\int_{t_2}^{\theta_3}\dots\int_{t_{M-1}}^{\theta_M} t_M e^{-\beta t_M} {dt_M}
\end{eqnarray*}}
{\color{blue}{Theorem statements}
\begin{itemize}
    
    \item A threshold policy  is always  in Best response  (Can we get uniqueness) 
    
    \item If uniqueness is done,  it is equivalent to work with threshold policies,  (try)
    
    \item Else, it is sufficient to work with threshold policies.. 
    
    \item  Monotonicity of thresholds 
    
    \item Define the Reduced game 
    
\end{itemize}
 
 We obtained the results of the reduced game through numerical simulations and discussed some interesting observations. 
}
}

 \end{document}